\documentclass{amsart}
\usepackage{amsxtra}
\theoremstyle{plain}

\newcommand{\cleqn}{\setcounter{equation}{0}}
\newcommand{\clth}{\setcounter{theorem}{0}}
\newcommand {\sectionnew}[1]{\section{#1}\cleqn\clth}
\newcommand{\nn}{\hfill\nonumber}
\newtheorem{theorem}{Theorem}[section]
\newtheorem{lemma}[theorem]{Lemma}
\newtheorem{definition-theorem}[theorem]{Definition-Theorem}
\newtheorem{proposition}[theorem]{Proposition}
\newtheorem{corollary}[theorem]{Corollary}
\newtheorem{definition}[theorem]{Definition}
\newtheorem{example}[theorem]{Example}
\newtheorem{remark}[theorem]{Remark}
\newtheorem{notation}[theorem]{Notation}
\newcommand \bth[1] { \begin{theorem}\label{t#1} }
\newcommand \ble[1] { \begin{lemma}\label{l#1} }

\newcommand \bpr[1] { \begin{proposition}\label{p#1} }
\newcommand \bco[1] { \begin{corollary}\label{c#1} }
\newcommand \bde[1] { \begin{definition}\label{d#1}\rm }
\newcommand \bex[1] { \begin{example}\label{e#1}\rm }
\newcommand \bre[1] { \begin{remark}\label{r#1}\rm }

\newcommand \bnota[1] { \begin{notation}\label{n#1}\rm }
\newcommand {\eth} { \end{theorem} }
\newcommand {\ele} { \end{lemma} }

\newcommand {\epr} { \end{proposition} }
\newcommand {\eco} { \end{corollary} }
\newcommand {\ede} { \end{definition} }
\newcommand {\eex} { \end{example} }
\newcommand {\ere} { \end{remark} }

\newcommand {\enota} { \end{notation} }
\newcommand \thref[1]{Theorem \ref{t#1}}
\newcommand \leref[1]{Lemma \ref{l#1}}
\newcommand \prref[1]{Proposition \ref{p#1}}

\newcommand \deref[1]{Definition \ref{d#1}}

\newcommand \reref[1]{Remark \ref{r#1}}
\newcommand \lb[1]{\label{#1}}

\def \Cset {{\mathbb C}}

\def \tP   {{\wt{P}}}

\def \tp {{\tilde{p}}}

\def \PG    {P(\Ga_1, \Ga_2)}
\def \A  {{\mathcal{A}}}           

\def \C  {{\mathcal{C}}}

\def \Z  {{\mathcal{Z}}}
\def \SS {{\mathcal{S}}}
\def \L  {{\mathcal{L}}}
\def \bD {{\bf D}}
\def \Da {{A_1}}
\def \Ra {{A_2}}
\def \De {\Delta}   

\def \al {\alpha}
\def \be {\beta}

\def \Ga {\Gamma}

\def \ra  {\rightarrow}           

\def \til {\tilde}
\def \hra {\hookrightarrow}
\def \sub {\subset}

\def \ol {\overline}
\def \wt {\widetilde}

\def \hs {\hspace{.2in}}

\def \Id { {\mathrm{Id}} }
\def \Stab { {\mathrm{Stab}} }

\def \BD   { {\mathrm{BD}} }
\def \MStab { \Stab_M }
\def \UStab { \Stab_U }
\def \newnew {{\mathrm{new}}}
\def \new { {(w_1, w_2)} }
\def \RAnew {{R^{\newnew}_{\A}(w_1, w_2) }}
\def \RCnew {{R^{\newnew}_{\C}(w_1, w_2) }}
\def \Aanew { {A_{1}^{\newnew}(w_1, w_2)} }
\def \Abnew { {A_{2}^{\newnew}(w_1, w_2)} }
\def \Canew { {C_{1}^{\newnew}(w_1, w_2)} }
\def \Cbnew { {C_{2}^{\newnew}(w_1, w_2)} }
\def \Kvv  { {K(v_1, v_2)} }
\def \Lvv  { {L(v_1, v_2)} }
\def \Jvv  { {J(v_1, v_2)} }
\def \jvv  { {j(v_1, v_2)} }
\def \Zvv  { {Z(v_1, v_2)} }
\def \WW { {}^{A_1} \! W_1^{C_1} \times 
{}^{A_2} \! W_2^{C_2} }
\def \Aavv { A_1(v_1, v_2) }

\def \Cbvv { C_2(v_1, v_2) }

\def \Ad { {\mathrm{Ad}} }

\def \g  {\mathfrak{g}}   

\def \l {\mathfrak{l}}

\def \k {\mathrm{k}}

\DeclareMathOperator \End { {\mathrm{End}} }

\DeclareMathOperator \Dom { {\mathrm{Dom}} }

\begin{document}
\title[Double cosets in reductive algebraic groups]
{On a class of double cosets in reductive algebraic groups}
\author[Jiang-Hua Lu]{Jiang-Hua Lu}
\address{
Department of Mathematics   \\
The University of Hong Kong \\
Pokfulam Road               \\
Hong Kong}
\email{jhlu@maths.hku.hk}
\author[Milen Yakimov]{Milen Yakimov}
\address{
Department of Mathematics \\
UC Santa Barbara \\
Santa Barbara, CA 93106, U.S.A.}
\email{yakimov@math.ucsb.edu}
\date{}
\begin{abstract}
We study a class of double coset spaces 
$R_\A \backslash G_1 \times G_2 /R_\C$, 
where $G_1$ and $G_2$ are connected reductive algebraic groups, 
and $R_\A$ and $R_\C$ are certain spherical subgroups of 
$G_1 \times G_2$ obtained by ``identifying'' 
Levi factors of parabolic subgroups in $G_1$ and $G_2$. 
Such double cosets naturally appear in the  symplectic leaf
decompositions of Poisson homogeneous spaces of complex
reductive groups with the Belavin--Drinfeld Poisson 
structures. They also appear in orbit decompositions of 
the De Concini--Procesi compactifications of semi-simple 
groups of adjoint type. We find explicit parametrizations of the 
double coset spaces and describe the double cosets as homogeneous 
spaces of $R_\A \times R_\C.$ We further show that all such double 
cosets give rise to set-theoretical solutions to the quantum 
Yang--Baxter equation on unipotent algebraic groups.
\end{abstract}
\maketitle
\sectionnew{The setup}\lb{setup}
\subsection{The setup}
\lb{adm}
 
Let $G_1$ and $G_2$ be two connected reductive algebraic groups 
over an algebraically closed base field $\k.$ 
For $i = 1, 2$, we will fix a maximal torus $H_i$ in $G_i$
and a choice $\Delta_{i}^{+}$ of positive roots in the set
$\Delta_i$ of all roots for $G_i$ with respect to $H_i$. For each 
$\alpha \in \Delta_{i}$, we will use $U_{i}^{\alpha}$  to denote the 
one-parameter unipotent subgroup of $G_i$ determined by $\alpha$. 
Let $\Gamma_i$ be the set of simple roots in $\Delta_{i}^{+}$. 
For a subset $A_i$ of $\Gamma_i$, let $P_{A_i}$ be the standard parabolic 
subgroup of $G_i$ containing the Borel subgroup $B_i=H_iU_i$, where
$U_i = \prod_{\alpha \in \Delta_{i}^{+}} U_{i}^{\alpha}$.
Let $M_{A_i}$ 
and $U_{A_i}$ be respectively the Levi factor of $P_{A_i}$ containing
$H_i$ and the unipotent radical of $P_{A_i}$. Then  
$P_{A_i}= M_{A_i} U_{A_i}$ and $M_{A_i}$ and $U_{A_i}$ 
intersect trivially. Denote by $M'_{A_i}$ the derived subgroups
of the Levi factors $M_{A_i}.$ We will also use
$\De_{A_i}$ to denote the set of roots in $\Delta_i$ that are 
in the linear span of $A_i$. 
Set $\De_{A_i}^+ = \De_{A_i}\cap \De^+.$
  
\bde{a-graph}
Identify $\Ga_1$ and $\Ga_2$ with the Dynkin diagrams of 
$G_1$ and $G_2$ respectively. 
By a {\it partial isometry} from $\Ga_1$ to $\Ga_2$ we mean 
an isometry from a subdiagram  of $\Ga_1$ to a subdiagram of $\Ga_2$. 
We will denote by $P(\Ga_1, \Ga_2)$ the set of all
partial isometries from $\Ga_1$ to $\Ga_2$. 
If $a \in \PG$, and if $\Da$ and $\Ra$ are the domain and the
range of $a$ respectively,  
we define a {\it generalized $a$-graph} to be an abstract
(not necessarily algebraic) 
subgroup $K$ of $M_{\Da} \times M_{\Ra}$ such that 
$K \cap (U_1^\alpha \times U_2^{a(\alpha)}) \subset 
U_1^\alpha \times U_2^{a(\alpha)}$
is the graph of a group isomorphism from 
$U_1^\alpha$ to $U_2^{a(\alpha)}$
for every $\alpha \in \De_{A_1}$. 
 
By an {\it admissible pair} for $G_1 \times G_2$
we mean a pair $(a, K)$, where
$a \in \PG$ and $K$ is a generalized $a$-graph. If
$\A=(a, K)$ is an admissible pair, we define the subgroup
$R_{\A}$ of $P_{\Da} \times P_{\Ra}$ by
\begin{equation}
\lb{R}
R_{\A} = K(U_{\Da} \times U_{\Ra}).
\end{equation}
\ede

\subsection{An Example}\lb{exe} 
Assume that $G_1$ and $G_2$ are connected and
simply-connected. Fix a
partial isometry $a \in P(\Ga_1, \Ga_2)$ with domain $A_1$ and 
range $A_2.$ Then
the derived subgroups $M_{A_i}^{\prime}$ of the Levi factors 
$M_{A_i}$ are connected and simply-connected 
\cite[Corollary 5.4]{SpSt}, and $a$
can be lifted to a group isomorphism 
$\theta_a \colon M_{A_1}^{\prime} \ra M_{A_2}^{\prime}.$ 
Let 
\[
{\rm Graph}(\theta_a)=
\{(m_1, \theta_a(m_1)) \mid m_1 \in M_{A_1}^{\prime}\}
\subset M_{A_1}^{\prime} \times M_{A_2}^{\prime}
\]
be the graph of $\theta_a$, and let $Z_{A_i}$ be the 
center of $M_{A_i}$
for $i = 1, 2$. Then,
any abstract subgroup $X$ of $Z_{A_1} \times Z_{A_2}$ 
gives rise to 
the group 
$K = X {\rm Graph}(\theta_a)$ 
which can be easily seen to be a generalized $a$-graph
for $G_1 \times G_2,$ and for the corresponding 
admissible pair $\A=(a, K)$, we have 
\[
R_\A= X {\rm Graph}(\theta_a) \, 
(U_{A_1} \times U_{A_2}).
\]
As a direct consequence of \leref{quintuple} in $\S$\ref{prp}, 
one sees that, in the case when $G_1$ and $G_2$ are connected and 
simply-connected, all admissible pairs $\A$ and groups $R_\A$
are of the above type. 
 
\subsection{Main problem}\lb{mainpr}
Fix two admissible pairs
$\A=(a, K)$ and $\C=(c, L)$.
In this paper, we obtain the solutions to the following problems:

{\em{1.) Find an explicit parametrization of all
$(R_{\A}, R_{\C})$ double cosets in $G_1 \times G_2.$

2.) Describe explicitly each $(R_{\A}, R_{\C})$ double coset
as a homogeneous space of $R_{\A} \times R_{\C}$.

3.) Find a closed formula for the dimension of each 
$(R_{\A}, R_{\C})$ double coset. 
}}

\bre{negative}
For an admissible pair $\A = (a, K)$, we will set
\begin{equation}
\lb{RAnegative}
R_{\A}^{-} = K (U_{A_1} \times U_{A_2}^{-}),
\end{equation}
where $U_{A_2}^{-}$ is the subgroup of $G_2$ generated by the 
one-parameter subgroups $U_{2}^{\alpha}$ for 
$\alpha \in - \De_{A_2}^+$. It is easy to see that solutions to the problems
in $\S$\ref{mainpr} lead to solutions to the corresponding
problems for the $(R_{\A}^{-}, R_{\C}^{-})$ double cosets in $G_1 \times G_2$.
\ere

\subsection{Motivation}\lb{motivation} 
Let $G$ be a  complex semi-simple Lie group of adjoint type.
The wonderful compactification $\ol{G}$ of $ G$, constructed by 
De Concini and Procesi in \cite{DP}, is a smooth projective 
$(G \times G)$-variety with finitely many $G \times G$ orbits. 
It was proved in \cite{DP} that every $G \times G$ orbit 
on $\ol{G}$ is isomorphic to 
$(G \times G)/R_{\C}^-$ for an
admissible pair $\C=(c, L)$ for $G \times G,$
where $c$ is the restriction of the identity automorphism 
of the Dynkin diagram of $G$ to some subdiagram.
Let $\g$ be the Lie algebra of $G$. De Concini
and Procesi further defined an embedding
of the wonderful compactification of $G$ in 
the projective variety of all Lie subalgebras
of $\g \oplus \g.$ Its image is the closure 
of the $G \times G$ orbit of the diagonal subalgebra
of $\g \oplus \g$ which clearly lies 
in the variety $\L_\g$ of all Lagrangian 
subalgebras of $\g \oplus \g$. Here, 
a subalgebra $\l$ of $\g \oplus \g$ is said to be
Lagrangian if $\dim \l = \dim \g$ and if $\l$ is
isotropic with respect to the symmetric 
bilinear form on $\g \oplus \g$ given by
\[
\langle x_1 + y_1, x_2 + y_2 \rangle = 
\ll x_1, x_2 \gg - \ll y_1, y_2 \gg, 
\hs x_1, y_1, x_2, y_2 \in \g,
\]
where  $\ll \, , \, \gg$ denotes the Killing form of $\g$.
The variety $\L_\g$ was recently studied in \cite{LuE2}
in connection with Drinfeld's classification \cite{D}
of Poisson
homogeneous spaces of Poisson Lie groups.
It was further shown in \cite{LuE2} that
all $G\times G$ orbits in $\L_\g$ (under the adjoint action)
are again of the 
type $(G \times G)/R_{\C}^{-}$ for some admissible pairs 
$\C=(c, L)$ for $G \times G,$ but this time the partial isometries
$c$ can be arbitrary. Thus, solutions to the double coset problems
in $\S$\ref{mainpr} can be used to classify 
$R_{\A}^{-}$ orbits in $\L_\g$ for any admissible pair $\A$.

Our motivation for studying $R_{\A}^{-}$ 
orbits in $\L_\g$ comes from
the theory of Poisson Lie groups.  In \cite{BD},
Belavin and Drinfeld classified all quasi-triangular Poisson
structures on $G$. Denote such a Poisson structure on $G$ by $\pi_{\BD}$.
The Poisson Lie groups $(G, \pi_{\BD})$ have received a lot of attention 
since the appearance of \cite{BD}. They were explicitly quantized
by Etingof--Schedler--Schiffmann in \cite{ESS} and their symplectic 
leaves were studied by the second author in \cite{Y}. 
One important problem regarding 
the Poisson Lie groups $(G, \pi_{\BD})$ is the study of their
Poisson homogeneous spaces \cite{D}. The variety $\L_\g$ of
Lagrangian subalgebras of $\g \oplus \g$ serves as a
``master" variety of Poisson homogeneous spaces of $(G, \pi_{\BD})$.
Indeed, identify $G$ with the diagonal
subgroup $G_\Delta$ of $G \times G$. It is shown in
\cite{LuE1} that $\L_\g$ carries a natural Poisson structure
$\Pi_{\BD}$, such that 
every $G_\Delta$ orbit in
$\L_\g$, when equipped with $\Pi_{\BD}$, is a Poisson homogeneous space of 
$(G, \pi_{\BD})$. To identify these homogeneous spaces of 
$G$, one needs to classify $G_\Delta$ orbits in
$\L_\g$, and to study their symplectic leaf decompositions with respect to 
$\Pi_{\BD}$, one
needs to classify $N(G_{\BD}^{*})$ orbits in $\L_\g$, where
 $G_{\BD}^{*}$ is the dual Poisson Lie group of
$(G, \pi_{\BD})$ inside $G \times G$, and $N(G_{\BD}^{*})$
is the normalizer subgroup of $G_{\BD}^{*}$ in $G \times G$. 
Both $G_{\BD}^{*}$ and $N(G_{\BD}^{*})$ are of the type $R_{\A}^{-}$,
where the partial isometry $a$ in $\A=(a, K)$ is the Belavin-Drinfeld 
triple defined in \cite{BD}. 
In the forthcoming paper \cite{LuY}, we will show that 
intersections of $G_\Delta$ and $N(G_{\BD}^{*})$ orbits in
$\L_\g$ are all regular Poisson subvarieties of
$(\L_\g, \Pi_{\BD})$. Solutions to the problems in $\S$\ref{mainpr} 
will be used in \cite{LuY}
to classify such orbits and to compute the ranks of the
Poisson structures $\Pi_{\BD}$. Based on the results of this paper,
in \cite{LuY} we will also treat the more general class of Poisson 
structures on $\L_\g,$ obtained from arbitrary lagrangian splittings
of $\g \oplus \g.$ The latter were classified by Delorme in 
\cite{De}. 

For the so-called standard Poisson structure on $G$, we have
$N(G_{\BD}^{*}) = B \times B^-$, where $(B, B^-)$ is 
a pair of opposite Borel subgroups of 
$G$. In \cite{Sp}, Springer classified the $B \times B^-$ orbits
in the wonderful compactification $\overline{G}$ of $G$ and 
studied the intersection cohomology of the $B \times B^-$ orbit closures. 
A classification of $G_\Delta$ orbits in $\overline{G}$ was
given by Lusztig in \cite{L2}. It was used by He \cite{H}
to study the closure of the unipotent variety
in $\overline{G}$. 
Our \thref{main1}  generalizes both Springer's and Lusztig's 
results.
We also point out that intersections of $G_\Delta$ and $B \times B^-$
orbits in the closed $G \times G$ orbits in $\L_\g$ are closely
related to the double Bruhat cells in $G$ and intersections of dual
Schubert cells on the flag variety of $G$ 
(see \cite{BGY, LuE2, FZ1}). It would be
very interesting to understand intersections of arbitrary 
$G_\Delta$ and $N(G_{\BD}^{*})$ orbits in $\L_\g$ in
the framework of the theory of cluster algebras
of Fomin and Zelevinsky \cite{FZ2}, and to find toric charts
on the symplectic leaves of $(\L_\g, \Pi_{\BD})$ using the methods of
Kogan and Zelevinsky \cite{KZ}.  

Our approach to the double coset space 
$R_\A \backslash G_1 \times G_2 /R_\C$  
relies on an inductive argument that relates an 
$(R_\A, R_\C)$ double coset in $G_1 \times G_2$ to 
cosets of similar type in the product $M_1 \times M_2$
of Levi subgroups  $M_1 \subset G_1$ and $M_2 \subset G_2.$ 
Such an argument first appeared in the work of the second author
\cite{Y} on the symplectic leaf decomposition of the
Poisson Lie group $(G, \pi_{\BD})$. It dealt with 
the case $G_1 = G_2 = G$ and $\A=\C$. Recently 
the first author and Evens \cite{LuE2} 
extended the methods of \cite{Y} to the case $G_1 = G_2 = G$ 
and $R_\A = G_\Delta.$ The advantage of the 
approach in this paper to the one in \cite{Y} 
is that in this paper we find an explicit relation between 
the double cosets in $G_1 \times G_2$ and the double cosets 
in $M_1 \times M_2$,
while \cite{Y} deals with various complicated 
Zariski open subsets of the   double 
cosets in $G_1 \times G_2$. In particular, the iteration of the inductive
procedure in \cite{Y} is only helpful for 
studying Zariski open subsets of the symplectic leaves
of $\pi_{\BD}$,
while we expect that the results of this paper
will be useful in understanding globally all the symplectic
leaves of $(\L_\g, \Pi_{BD}).$ 

The roots of the methods used in  this paper
can be traced back to the inductive procedure
of Duflo \cite{Du} and Moeglin--Rentschler \cite{MR}
for describing 
the primitive spectrum of the universal enveloping algebra
of a Lie algebra which is in general neither solvable nor 
semi-simple. The reason for this relation
is that the Poisson structures
$\Pi_{BD}$ on $\L_\g$ vanish at many points
and the corresponding linearizations give rise
to Lie algebras that are rarely semi-simple
or solvable. Due to the Kirillov--Kostant 
orbit method, one expects a close relation between
the primitive spectrum of the quantized algebras
of coordinate rings of various affine subvarieties
of $\L_\g$ (which are in fact deformations of the 
universal enveloping algebras of the linearizations)
and the symplectic leaves of such affine Poisson 
subvarieties of $\L_\g.$ From this point of view, 
our method can be considered 
as a nonlinear quasiclassical analog of \cite{Du, MR}. 

\subsection{Organization of paper}\lb{org}
The statements of the main results in this paper are given in
$\S$\ref{results}. The two sections, $\S$\ref{prp} and $\S$\ref{facts}, 
serve as preparations for the proofs of the main 
results which are given in
$\S$\ref{prmain1} and $\S$\ref{prmain2}. 
In $\S$\ref{YBE}, we describe some solutions to the 
set-theoretical Yang--Baxter equation that arise in our solution 
to 2) of the main problems in $\S$\ref{mainpr}.
In $\S$\ref{spherical}, we classify 
$(R_\A, P)$ double cosets in 
$G_1 \times G_2$  for any
admissible pair $\A$ and any 
standard parabolic subgroup $P$ of $G_1 \times G_2$,
and we show in particular that $R_\A$ is a spherical
subgroup of $G_1 \times G_2.$ 

\subsection{Notation}\lb{notation}
Throughout this paper, the letter $e$ always denotes the 
identity element of any group. We will denote by $\Ad$ 
the conjugation action of a group $G$ on itself,
given by $\Ad_g(h) = ghg^{-1}$ for $g, h \in G.$ 
If $M_1$ and $M_2$ are two sets and if $K \subset M_1 \times M_2$, 
we can think of $K$ as a correspondence between $M_1$ and $M_2$. 
If $M_3$ is another set and if $L \subset M_2 \times M_3$, then one
defines the composition of $L$ and $K$ to be
\begin{equation}
\lb{composition}
L \circ K = \{(m_1, m_3)\in M_1 \times M_3 \mid \exists m_2 \in M_2
\; \; \mbox{such 
that} \; \; (m_1, m_2) \in K, \; (m_2, m_3) \in L\}.
\end{equation}

\subsection{Acknowledgements}\lb{ack} This work was started
while the authors were visiting the Erwin Schr\"odinger institute for 
Mathematical
Physics, Vienna in August 2003. We would like to thank the organizers 
A. Alekseev and T. Ratiu of the program on Moment Maps 
for the invitation to participate. During the
course of the work, the first author enjoyed the hospitality of IHES
and UCSB and the second author of the University of Hong Kong. 
We are grateful to these institutions for the warm 
hospitality. We would also 
like to thank S. Evens, K. Goodearl, G. Karaali and A. Moy 
for helpful discussions. We are indebted to the referee
for carefully reading the paper, pointing out many
typos, and suggesting several improvements in the exposition.
The research of the first author was partially supported by 
(USA)NSF grant DMS-0105195 and HKRGC grant 701603.
The research of the second author was partially supported by 
NSF grant DMS-0406057 and
a UCSB junior faculty research incentive grant.   

\sectionnew{Statements of the main results}
\lb{results}
 
The main parameters that come into our classification of
$(R_{\A}, R_{\C})$ double cosets in
$G_1 \times G_2$ are 
certain elements in the Weyl groups of $G_1$ and $G_2$ and
certain
``twisted conjugacy classes"
in Levi subgroups of $G_1$. 
 
\subsection{Minimal length representatives}
\lb{minrep}
 
For $i = 1, 2$, let $W_i$ be the Weyl group of $\Ga_i$. For 
a subset $D_i$ of $\Ga_i$, we will use $W_{D_i}$ to denote  the
subgroup of $W_i$ generated by $D_i$. Let $W_i^{D_i} \subset W_i$ and
${}^{D_i}\!W_i \subset W_i$ be respectively the set of minimal length 
representatives
of the cosets from $W_i/W_{D_i}$ and $W_{D_i} \backslash W_i$.
It is well-known that
\[
w_i \in W_i^{D_i} \; 
\hs {\mbox{if and only if}} \hs
w_i(D_i) \subset \Delta_{i}^{+} \; 
\]
and
\[
w_i \in {}^{D_i}\!W_i  \hs {\mbox{if and only if}} \hs
w_{i}^{-1}(D_i) \subset \Delta_{i}^{+}.
\]
 
\subsection{Twisted conjugacy classes}
\lb{tc}
 
\bde{tc}
Let $D_1$ be a subset of $\Ga_1$, and let $d$ be a 
partial isometry  from $\Ga_1$ to $\Ga_1$ with $D_1$ 
as both its domain and its range. If 
$J \subset M_{D_1} \times M_{D_1}$ is a generalized $d$-graph, we let
$J$ act on $M_{D_1}$ from the left by
\[
J \times M_{D_1} \ra M_{D_1} \colon ((l, n), m) \mapsto 
l m n^{-1}, \hs (l, n) \in J, \; m \in M_{D_1}.
\]
By a {\it $d$-twisted conjugacy class} in $M_{D_1}$ we mean
a $J$ orbit in $M_{D_1}$  for
any generalized $d$-graph $J$. If $J$ is the graph of
a group automorphism 
$j \colon M_{D_1} \to M_{D_1}$, the action of
$J$ on $M_{D_1}$ becomes the usual {\it $j$-twisted 
conjugation action} of
$M_{D_1}$ on itself:
\[
M_{D_1} \times M_{D_1} \ra M_{D_1} \colon (l, m) \mapsto 
l m j(l)^{-1}, \hs l, m \in M_{D_1}
\]
whose orbits will be referred to as
{\it $j$-twisted conjugacy classes}.
\ede
   
\subsection{Classification of double cosets}
\lb{classification}

Let $\A = (a, K)$ and $\C = (c, L)$ be two admissible pairs, and let
$A_1$ and $A_2$ be the domain and the 
range of $a$ and $C_1$ and $C_2$
the domain and the  range of $c$ respectively. We now state our 
classification
of $(R_\A, R_\C)$ double cosets in $G_1 \times G_2$. In doing so,
we will need to refer to several lemmas that will be proved in 
$\S$\ref{prp}-$\S$\ref{prmain2}.
   
For $v_1 \in W_1^{C_1}$ and $ v_2 \in {}^{A_2} \! W_2$, consider the 
sets
\begin{align}
\lb{A1stab}
A_1(v_1, v_2) =\{ \al \in A_1 \mid
&(v_1 c^{-1} v_2^{-1} a)^n \al \; \mbox{is defined}
\\ 
\nn
&\mbox{and is in} \; A_1 \; \mbox{for} \; n=0, 1, \ldots\}
\\
\lb{C2stab}
C_2(v_1, v_2) =\{ \beta \in C_2 \mid
&(v_{2}^{-1} a v_1 c^{-1})^n \beta \; \mbox{is defined}
\\ 
\nn
&\mbox{and is in} \; C_2 \;\mbox{for} \; n=0, 1, \ldots\}
\end{align}
Then
$A_1(v_1, v_2)$ is the largest subset of $A_1$ that is invariant under
$v_1 c^{-1} v_2^{-1} a$, and $C_2(v_1, v_2)$
is the largest subset of $C_2$ that is invariant under
$v_{2}^{-1}a v_1 c^{-1}$. Note that
\begin{equation}
\lb{st_isom}
v_{2}^{-1} a \; \;  \mbox{and} \; \;  
cv_{1}^{-1} \colon A_1(v_1, v_2) \ra 
C_2(v_1, v_2)
\end{equation}
are both partial isometries from $\Ga_1$ to $\Ga_2$. 
Fix a representative $\dot{v}_i$ in the normalizer of $H_i$ in $G_i$, 
and define 
\begin{align}
\lb{Kvv}
\Kvv &= \left(M_{A_1(v_1, v_2) } \times 
M_{C_2(v_1, v_2)}\right) \cap \Ad_{(e,\dot{v}_2)}^{-1} K\\
\lb{Lvv}
\Lvv &= \left(M_{A_1(v_1, v_2)} 
\times M_{C_2(v_1, v_2)}\right)
\cap \Ad_{(\dot{v}_1, e)} L.
\end{align} 
It is easy to see that $\Kvv$ is a generalized $(v_{2}^{-1}a)$-graph, 
while $\Lvv$ is a generalized $(cv_{1}^{-1})$-graph for the partial 
isometries in \eqref{st_isom}.
Define 
\begin{align}
\lb{Jvv}
\Jvv = \{&(m, m^\prime) \in M_{A_1(v_1, v_2)} \times
M_{A_1(v_1, v_2)} \mid \exists n \in M_{C_2(v_1, v_2)} \;
\\
\nn
&\mbox{such that} \;
(m, n) \in \Kvv \; \mbox{and} \; 
(m^\prime, n) \in \Lvv\}.
\end{align}
Then $\Jvv$ is a generalized
$(v_1 c^{-1} v_2^{-1} a)$-graph for the partial isometry 
\[
v_1 c^{-1} v_2^{-1} a \colon A_1(v_1, v_2) \ra A_1(v_1, v_2).
\]
This follows from \leref{composition} in $\S$\ref{prp}
and the facts that if
$\sigma$ is the map
\[
\sigma \colon M_{A_1(v_1, v_2)} \times M_{C_2(v_1, v_2)} \ra
M_{C_2(v_1, v_2)} \times M_{A_1(v_1, v_2)} \colon 
(m, n) \mapsto (n, m),
\]
then $\sigma(\Lvv)$ is a generalized $(v_1 c^{-1})$-graph, and
$\Jvv$ is the composition
\begin{equation}
\lb{compJ}
\Jvv= \sigma(\Lvv) \circ \Kvv
\end{equation}
(of group correspondences, see $\S$\ref{notation}).
Let $\Jvv$ act on $M_{A_1(v_1, v_2)}$ from the left by
\begin{equation}
\lb{Q1-M1}
(l_1, n_1) \cdot m_1 := l_1 m_1 n_{1}^{-1}, \hspace{.2in}
(l_1, n_1) \in \Jvv, \; m_1 \in M_{A_1(v_1, v_2)}.
\end{equation}
By \deref{tc}, $\Jvv$ orbits in $M_{A_1(v_1, v_2)}$
are $(v_1 c^{-1} v_2^{-1} a)$-twisted conjugacy classes 
in $M_{A_1(v_1, v_2)}$.
Let $Z_{C_2(v_1, v_2)}$ be the center of 
$M_{C_2(v_1, v_2)}$, let  
$\eta_2 \colon G_1 \times G_2 \to G_2$ be the projection, and let
\[
\Zvv =Z_{C_2(v_1, v_2)}/ \left(Z_{C_2(v_1, v_2)} \cap 
 v_{2}^{-1}(\eta_2(K))\right)
\left(Z_{C_2(v_1, v_2)} \cap \eta_2(L)\right).
\]
For each $s \in \Zvv$, 
fix a representative $\dot{s}$ of $s$ in $Z_{C_2(v_1, v_2)}$.
For $(g_1, g_2) \in G_1 \times G_2$,
we will use $[g_1, g_2]$ to denote the double coset
$R_\A(g_1, g_2)R_\C$ in $G_1 \times G_2$.
We can now state our main theorem on the 
classification of $(R_\A, R_\C)$ 
double cosets in $G_1 \times G_2$. 
  
\bth{main1} Let  $\A=(a, K)$ and $\C=(c, L)$ be 
any two admissible
pairs for $G_1 \times G_2$.
Then every $(R_\A, R_\C)$ double coset of $G_1 \times G_2$ is of the form
\[
[m_1 \dot{v}_1, \dot{v}_2 \dot{s}] \; \mbox{for some} \;
v_1 \in W_1^{C_1}, v_2 \in {}^{A_2} \! W_2, \;
m_1 \in M_{A_1(v_1, v_2)}, 
{s} \in \Zvv.
\]
Two such double cosets $[m_1 \dot{v}_1, \dot{v}_2 \dot{s}]$ and 
$[m^\prime_1 \dot{v}'_1, \dot{v}'_2 \dot{s}^\prime]$ coincide 
if and only if $v'_i=v_i$ for $i = 1, 2$, $s = s^\prime$, and 
$m_1$ and $m^\prime_1$ are in the same $\Jvv$ orbit in
$M_{A_1(v_1, v_2)}$.
\eth

\bex{full} If $\eta_2(K) = M_{A_2}$ or $\eta_2(L) = M_{C_2}$, then $\Zvv=\{e\}$ 
for all choices of $(v_1, v_2).$ In particular, 
every $(R_\A, R_\C)$ double coset of $G_1 \times G_2$ is of the form
\[
[m_1 \dot{v}_1, \dot{v}_2] \; \mbox{for some} \;
v_1 \in W_1^{C_1}, v_2 \in {}^{A_2} \! W_2, \;
m_1 \in M_{A_1(v_1, v_2)}.
\]
Two cosets of this type $[m_1 \dot{v}_1, \dot{v}_2]$ and 
$[m^\prime_1 \dot{v}'_1, \dot{v}'_2]$ coincide 
if and only if $v'_i=v_i$ for $i = 1, 2$, and 
$m_1$ and $m^\prime_1$ are in the same $\Jvv$ orbit in
$M_{A_1(v_1, v_2)}$.
\eex

\bex{graph}
If $K$ and $L$ are 
respectively the graphs of group isomorphisms
\[
\theta_a \colon M_{A_1} \to M_{A_2} \hspace{.2in} {\rm and} \hspace{.2in}
\theta_c \colon M_{C_1} \to M_{C_2},
\]
the group $\Jvv$ is the graph of the group automorphism
\[
\jvv:=\Ad_{\dot{v}_1} \theta_{c}^{-1} 
\Ad_{\dot{v}_{2}}^{-1} \theta_a \colon 
M_{A_1(v_1, v_2)} \to M_{A_1(v_1, v_2)}.
\]
Then every $(R_\A, R_\C)$ orbit in 
$G_1 \times G_2$ is of the form $[m_1 \dot{v}_1, \dot{v}_2]$ for 
a unique pair $(v_1, v_2) \in W_1^{C_1} \times {}^{A_2} \! W_2$
and some $m_1 \in M_{A_1(v_1, v_2)}$. Two such cosets
$[m_1 \dot{v}_1, \dot{v}_2]$ and 
$[m^\prime_1 \dot{v}_1, \dot{v}_2]$
coincide if and only if $m_1$ and $m^\prime_1$ are in the same
$\jvv$-twisted 
conjugacy class (see $\S$\ref{tc}).
\eex
\subsection{Structure of double cosets}\lb{stabilizer}
Let $\A = (a, K)$ and
$\C = (c, L)$ be two admissible pairs for $G_1 \times G_2$. 
For $q = (g_1, g_2) \in G_1 \times G_2$, 
set
\[
\Stab(q) = \Ad_{(e, g_2)}^{-1}(R_\A) \cap \Ad_{(g_1, e)} (R_\C).
\]
Then the double coset $R_\A q R_\C$ in $G_1 \times G_2,$ 
considered as an $R_{\A} \times R_{\C}$-space under the action
$(r_{\A}, r_{\C}) \cdot q'= r_{\A}q' r_{\C}^{-1},$
is given by
\[
R_\A q R_\C = (R_\A \times R_\C) /
\{( \Ad_{(e, g_2)}(r), \Ad_{(g_1, e)}^{-1}(r))\mid
r \in \Stab(q)\}.
\]
For
$q = (m_1\dot{v}_1, \dot{v}_2 s_2)$ as in 
\thref{main1}, where $m_1 \in M_{\Aavv}, 
v_1 \in W_{1}^{C_1}, v_2 \in {}^{A_2}\!W_2$, and 
$s_2 \in Z_{\Cbvv}$, we will now give an explicit description of $\Stab(q)$.
Let $\pi_{A_1} \colon P_{A_1} \to M_{A_1}$ and 
$\pi_{C_2} \colon P_{C_2} \to M_{C_2}$ be, respectively, the projections 
with respect to the decompositions $P_{A_1} = M_{A_1}N_{A_1}$
and $P_{C_2} = M_{C_2}N_{C_2}$. By
2) of \leref{property}  in $\S$\ref{prp} (by taking $D_1 = \emptyset$
therein),
the group $K$ gives rise to a
group isomorphism 
$\theta_a \colon U_1 \cap M_{A_1} \to U_2 \cap M_{A_2}$
such that $\theta_a(U_{1}^{\alpha}) = U_{2}^{a(\alpha)}$ for every 
$\alpha \in \De_{A_1}$. Similarly, we have the group isomorphism
$\theta_c \colon U_1 \cap M_{C_1} \to U_2 \cap M_{C_2}$ 
induced from the group
$L$. We will prove that
\begin{equation}
\lb{st1}
\Stab(q) \subset
\left(P_{\Aavv} \cap v_1(P_{C_1})\right) \times 
\left(P_{\Cbvv} \cap v_2^{-1}(P_{A_2}) \right),
\end{equation}
and we will use the maps $\pi_{A_1}, \pi_{C_2}, \theta_a$,
and $\theta_c$ to describe $\Stab(q)$.

Note that $v_1 \in {}^{\Aavv}\!W_{1}^{C_1}$ and 
$v_2 \in {}^{A_2}\!W_{2}^{\Cbvv}$ because 
\begin{equation}\lb{AC}
v_{1}^{-1}(\Aavv) = c^{-1}(\Cbvv) \subset C_1  \; \; \; 
\mbox{and}
\; \; \; v_2(\Cbvv) = a(\Aavv) \subset A_2.
\end{equation}
By \leref{int} in $\S$\ref{facts}, 
we have
\begin{align}
P_{\Aavv} \cap v_1(P_{C_1})& = M_{\Aavv}
 \left(U_{\Aavv} \cap v_1(U_1)\right),
\lb{st2}
\\
P_{\Cbvv} \cap v_2^{-1}(P_{A_2})& =
M_{\Cbvv} \left(U_{\Cbvv} \cap v_2^{-1}(U_2)\right).
\lb{st3}
\end{align}
Consider the group homomorphisms
\begin{align}
\lb{phi1}
\phi_1 &= \Ad_{m_1 \dot{v}_1} \theta_{c}^{-1} \pi_{C_2} \colon \;
U_2 \ra \Ad_{m_1}(U_1),\\
\lb{phi2}
\phi_2 &= \Ad_{\dot{v}_2s_2}^{-1} \theta_a \pi_{A_1} \colon \;U_1 \ra U_2.
\end{align}
In $\S$\ref{prmain2} we will show that 
\begin{align}
\lb{sigma1}
\sigma_1 &=\phi_1 \phi_2 \in 
\End(U_{\Aavv} \cap v_1(U_1)), \\
\lb{sigma2}
\sigma_2 &= \phi_2 \phi_1 \in
\End( U_{\Cbvv} \cap v_2^{-1}(U_2))
\end{align}
are well-defined and  that $\sigma_{1}^{k+1} \equiv e$ and 
$\sigma_{2}^{k+1}\equiv e$ for some integer $k \geq 1$,
where $e$ stands for the trivial endomorphism. 
Consider the subgroups 
\begin{align}
\lb{subgr1}
&U_1 \cap v_1(U_{C_1}) 
= U_{\Aavv} \cap v_1(U_{C_1}) 
\subset U_{\Aavv} \cap v_1(U_1), \\  
\lb{subgr2}
&U_2 \cap v_2^{-1}(U_{A_2})
= U_{\Cbvv} \cap v_2^{-1}(U_{A_2})
\subset U_{\Cbvv} \cap v_2^{-1}(U_2),
\end{align}
and define
\begin{align*}
\psi_1 & \colon \left(U_1 \cap v_1(U_{C_1})\right) 
\times \left(U_2 \cap v_2^{-1}(U_{A_2})\right) \ra 
U_{\Aavv} \cap v_1(U_1),\\
\psi_2 & \colon \left(U_1 \cap v_1(U_{C_1})\right) 
\times \left(U_2 \cap v_2^{-1}(U_{A_2})\right) \ra  
U_{\Cbvv} \cap v_2^{-1}(U_2),
\end{align*}
respectively by
\begin{align}
\lb{psi1}
\psi_1(n_1, n_2) &= \phi_1(n_2) \sigma_1(n_1 \phi_1(n_2)) 
\sigma_1^2(n_1 \phi_1(n_2)) \cdots \sigma_1^k(n_1 \phi_1(n_2)),\\
\lb{psi2}
\psi_2(n_1, n_2) &= \phi_2(n_1) \sigma_2(n_2 \phi_2(n_1)) 
\sigma_2^2(n_2 \phi_2(n_1)) \cdots \sigma_2^k(n_2 \phi_2(n_1)).
\end{align}

\bth{main2}
Let the notation be as in \thref{main1}. For
$q = (m_1\dot{v}_1, \dot{v}_2 s_2)$, where $m_1 \in M_{\Aavv}, 
v_1 \in W_{1}^{C_1}, v_2 \in {}^{A_2}\!W_2$, and 
$s_2 \in Z_{\Cbvv},$ we have the semi-direct product
decomposition
\[
\Stab(q) = \MStab(q) \UStab(q),
\]
where 
\begin{align*}
\MStab(q)& =\Kvv \cap \Ad_{(m_1, e)}\Lvv\\
&=\left(M_{\Aavv} \times M_{\Cbvv}\right)
\cap \Ad_{(e,\dot{v}_2 )}^{-1} K \cap \Ad_{(m_1\dot{v}_1, e)} L,\\
\UStab(q) &= \{ \left( n_1 \psi_1(n_1, n_2), \; n_2 \psi_2(n_1, n_2) \right) 
\mid 
(n_1, n_2) \in 
\left(U_1 \cap v_1(U_{C_1})\right) 
\times \left(U_2 \cap v_2^{-1}(U_{A_2})\right)\} \\
&\subset \left(U_{\Aavv} \cap v_1(U_1)\right) \times
\left(U_{\Cbvv} \cap v_2^{-1}(U_2)\right).
\end{align*}
In particular $\Stab(q)$ is the semi-direct product of
\begin{align*}
\Stab_M(q) &= \Stab(q) \cap \left( M_{\Aavv} \times M_{\Cbvv} \right)
\; \mbox{and} \\
\Stab_U(q) &= \Stab(q) \cap 
\left((U_{\Aavv} \cap v_1(U_1)) \times 
(U_{\Cbvv} \cap v_2^{-1}(U_2))\right),
\end{align*}
cf. \eqref{st1} and \eqref{st2}--\eqref{st3}.
\eth  

\bre{Jm}
For $m_1 \in M_{\Aavv}$, let $J_{m_1}(v_1, v_2)$ be the
stabilizer subgroup of $J(v_1, v_2)$ at $m_1$ for the
action of $\Jvv$ on $M_{\Aavv}$ given by 
\eqref{Q1-M1}, cf. \deref{tc}. We will show in 
\leref{ext} in $\S$\ref{prmain2} that $\MStab(q)$
is a central extension of $J_{m_1}(v_1, v_2)$ by
\[
(\ker \eta_1|_L) \cap (\Id \times v_2^{-1})(\ker \eta_1|_K)
\subset \{e\} \times Z_{\Cbvv},
\]
where, for $i = 1, 2$, 
$\eta_i \colon G_1 \times G_2 \to G_i$ is the
projection. 
\ere

\subsection{A dimension formula}\lb{dim}

\bth{dim} In the notation of \thref{main2}, if  
$K$ and $L$ are algebraic subgroups of $G_1 \times G_2$,
or if the base field is $\Cset$ and $K$ and $L$ are Lie subgroups
of $G_1 \times G_2$, then the double coset
$[m_1 \dot{v}_1, \dot{v}_2 \dot{s}]$ (which is clearly smooth 
subvariety or submanifold, respectively) has dimension equal to
\begin{align*}
  & l(v_1) + l(v_2)+ 
\dim P_{A_1} - \dim M_{A_1(v_1,v_2)} + 
\dim P_{C_2} - \dim Z_{C_2(v_1,v_2)}
\\
&+ \dim \eta_2
\left(KL \cap (v_1^{-1}Z_{A_1(v_1,v_2)}\times v_2 Z_{C_2(v_1,v_2)})
\right)
+ \dim (\Jvv \cdot m_1).
\end{align*}
Here $l(\cdot)$ denotes the length function on the 
Weyl groups $W_1$ and $W_2$, and $\Jvv \cdot m_1$ is the $\Jvv$ orbit
through $m_1$ for the $\Jvv$ action on $M_{\Aavv}$ 
given in \eqref{Q1-M1}.
\eth
 
Let us note that in the two cases in \thref{dim}, the double coset 
$[m_1 \dot{v}_1, \dot{v}_2 \dot{s}]$ is
a locally closed algebraic subset of $G_1 \times G_2$
(as an orbit of an algebraic group) 
or a submanifold of $G_1 \times G_2.$ It will be shown in $\S$\ref{prmain2}
that 
in these two cases $\Jvv$ is respectively an algebraic group or a Lie group. 
Note that 
$\eta_2\left(KL \cap (v_1^{-1} Z_{A_1(v_1,v_2)}\times v_2 Z_{C_2(v_1,v_2)}) 
\right)$ 
is also an algebraic group or a Lie group because it is a 
quotient of 
$KL \cap (v_1^{-1} Z_{A_1(v_1,v_2)}\times v_2 Z_{C_2(v_1,v_2)}).$
  
\sectionnew{Properties of the groups $R_\A$}
\lb{prp}
Let $a \in P(\Ga_1, \Ga_2)$ be a partial isometry from 
$\Ga_1$ to $\Ga_2$, and let
$A_1$ and $A_2$ be respectively the domain and the range of
$a$. In this section, we first give a 
characterization of those subgroups $K$ of $M_{A_1} \times M_{A_2}$ that
are generalized $a$-graphs. We will then prove some
properties of the groups $R_\A$ associated to an admissible pairs
$\A=(a, K)$. These properties are crucial in 
the proofs of the main results in this paper.

Recall that for $i = 1, 2$,  
$\eta_i \colon G_1 \times G_2 \to G_i$ denote the 
projections to the i'th factor. For a subset $D_i$ of $\Ga_i$, 
$M_{D_i}^{\prime}$ denotes the derived subgroup of $M_{D_i}.$ 
We will denote by $Z_{D_i}$ the center of $M_{D_i}$.
 
\bde{a-quintuples}
By an $a$-quintuple we mean a quintuple
$(K_1, X_1, K_2, X_2, \theta)$, where,
for $i = 1, 2$, $K_i$ is an abstract subgroup of $M_{A_i}$
containing $M_{A_i}^{\prime}$, $X_i$ is an abstract subgroup of
$K_i \cap Z_{A_i}$, and
$\theta \colon K_1/X_1 \ra K_2/X_2$ is a group isomorphism that
maps 
$U_{1}^{\alpha}$ to $U_{2}^{a(\alpha)}$ for each
$\alpha \in \De_{A_1}$, where $U_{1}^{\alpha}$ is identified with its 
image in $K_1/X_1$ and similarly for $U_{2}^{a(\alpha)}.$
\ede
 
\ble{quintuple}
Let $K \subset M_{A_1} \times M_{A_2}$ be a generalized
$a$-graph. Define
\begin{align*}
&K_1 = \eta_1(K)\subset M_{A_1}, \; \; \;
X_1 = \eta_1(\ker (\eta_2|_K)) = \{x_1 \in M_{A_1}
\mid (x_1, e) \in K\} \subset K_1\\
& K_2= \eta_2(K) \subset M_{A_2}, \; \; \;
X_2 = \eta_2(\ker (\eta_1|_K)) = \{x_2 \in
M_{A_2} \mid (e, x_2) \in K\}\subset K_2,
\end{align*}
and let
\begin{equation}
\lb{theta}
\theta \colon K_1/X_1 \ra K_2/X_2 \; \;
\mbox{where} \; \;  \theta(k_1X_1) = k_2X_2 \hs
{\mbox{if}} \; \; (k_1, k_2) \in K.
\end{equation}
Then $\theta$ is well-defined, and
$\Theta(K):=(K_1, X_1, K_2, X_2, \theta)$ is 
an $a$-quintuple. Moreover, $K$ can be expressed in terms of 
$\Theta(K)$
as
\begin{equation}
\lb{K}
K = \{(k_1, k_2) \in K_1 \times K_2 \; \mid \; 
\theta(k_1X_1) = k_2X_2\}.
\end{equation} 
The assignment $K \mapsto \Theta(K)$ is a 
one to one correspondence between
the set of generalized $a$-graphs and the set of
$a$-quintuples.
\ele
 
\begin{proof} Let $i=1,2.$ Since $K_i$ contains all one-parameter
unipotent subgroups $U^{\alpha}_i$ of $M_{A_i}$ and the latter generate
$M'_{A_i}$ as an abstract group, we have that $K_i \supset M'_{A_i}.$ 
Moreover $X_i$ is an abstract normal subgroup of $K_i$ 
which does not intersect any one-parameter unipotent subgroup 
$U_i^\alpha$ of $M_{A_i}$ because of the main condition for an 
$a$-graph, cf. \deref{a-graph}.  

First we show that $X_i \subset Z_{A_i}.$ We will make use of the 
following fact which can be found
e.g. in \cite[Corollary 29.5]{Hum}. 

{\em{$(*)$ Any simple algebraic group of adjoint type (over an algebraically 
closed field) is simple as an abstract group.}}

Assume that $X_i$ is not a subgroup  of $Z_{A_i}.$ 
Because of $K_i \supset M'_{A_i},$ we have that 
$M_{A_i}=K_iZ_{A_i}$. Since
$X_i$ is normal in $K_i,$ $X_iZ_{A_i}/Z_{A_i}$ is a nontrivial 
abstract normal subgroup of $M_{A_i}/Z_{A_i}.$ 
The group $M_{A_i}/Z_{A_i}$ is semi-simple and has trivial center, 
so it is a direct product of simple algebraic groups
of adjoint types. We can now apply $(*)$ to get
that $X_i Z_{A_i}/Z_{A_i}$ contains a simple factor of 
$M_{A_i}/Z_{A_i}.$ Therefore $(X_i Z_{A_i})'$ contains 
$M'_{D_i}$ for some nontrivial subset $D_i$ of $A_i$.
(Here $(.)'$ refers to the derived subgroup of an 
abstract group.) As a consequence 
$X_i \supset X'_i = (X_i Z_{A_i})' \supset M'_{D_i}$
which contradicts with the fact that $X_i$ does not intersect 
any one-parameter unipotent subgroup of $M_{A_i}.$ 
This shows that $X_i \subset Z_{A_i}.$  

It is easy to see that $\theta$ is a well-defined group isomorphism 
and that \eqref{K} holds. The definition of $K$ implies that
$\theta$ induces a group isomorphism from 
$U_{1}^{\alpha}$ to $U_{2}^{a(\alpha)}$ for each
$\alpha \in \De_{A_1}$. Thus $\Theta(K)$ is an $a$-quintuple.
It is straightforward to check that 
the map $K \mapsto \Theta(K)$ is bijective.
\end{proof}
 
\bnota{NAD}
For $i = 1, 2$ and for a subset $D_i \subset A_i$, let 
\[
P_{D_i}^{A_i} = P_{D_i} \cap M_{A_i}, \hs
U_{D_i}^{A_i} = U_{D_i} \cap M_{A_i}.
\]
For subsets $S_i \subset M_{A_i}$, we will set
\begin{eqnarray}
\lb{Kof}
K(S_1) &= \{m_2 \in M_{A_2} \mid \exists \, m_1 \in S_1 \; \; 
\mbox{such that} \; \; (m_1, m_2) \in K\}\\
K(S_2) &= \{m_1 \in M_{A_1} \mid \exists \, m_2 \in S_2 \; \; 
\mbox{such that} \; \; (m_1, m_2) \in K\}.
\end{eqnarray}
\enota
 
\ble{property}
For any subset $D_1$ of $A_1$, 
 
1) $K(M_{D_1}) \subset M_{a(D_1)}$, $K(M_{a(D_1)}) \subset
M_{D_1}$, and the intersection 
$(M_{D_1} \times M_{a(D_1)}) \cap K$ is generalized 
$a|_{D_1}$-graph;

2) the intersection $(U_{D_1}^{A_1} \times U_{a(D_1)}^{A_2}) \cap K$ 
is the graph of a group 
isomorphism
$\theta \colon U_{D_1}^{A_1} \to U_{a(D_1)}^{A_2}$;

3) we have the decompositions
\begin{eqnarray}
\lb{PK1}
(P_{D_1} \times M_{A_2}) \cap K &=&
(M_{A_1} \times P_{a(D_1)}) \cap K =
(P_{D_1} \times P_{a(D_1)} ) \cap K
\\
\lb{PK2}
&= &\left((M_{D_1} \times M_{a(D_1)})\cap K\right) 
\left((U_{D_1}^{A_1} \times U_{a(D_1)}^{A_2}) \cap K \right)\\
\lb{PR1}
(P_{D_1} \times P_{A_2}) \cap R_\A & = & 
(P_{A_1} \times P_{a(D_1)})\cap R_\A = 
(P_{D_1} \times P_{a(D_1)})\cap R_\A\\
\lb{PR2}
&=&
\left((M_{D_1} \times M_{a(D_1)})\cap R_\A\right) 
\left((U_{D_1} \times U_{a(D_1)}) \cap R_\A \right)\\
\lb{PR3}
&=&\left((P_{D_1} \times P_{a(D_1)} ) \cap K\right) 
\left(U_{A_1} \times U_{A_2} \right).
\end{eqnarray}
\ele

\begin{proof} We first prove that $K(H_1) \subset H_2$ (corresponding to
the special case when $D_1 = \emptyset$). Assume that 
$(h_1, m_2) \in (H_1 \times M_{A_2}) \cap K$. 
For $\beta \in \De_{A_2}$, let 
$\alpha = a^{-1}(\beta)$, and let 
$\theta_\alpha \colon U_{1}^{\alpha} \to U_{2}^{\beta}$ 
be the group isomorphism whose
graph is $(U_{1}^{\alpha}\times U_{2}^{\beta}) \cap K$. For
every $(u_1, u_2) \in 
(U_{1}^{\alpha}\times U_{2}^{\beta}) \cap K$, we have 
\[
(h_1 u_1 h_{1}^{-1}, \; m_{2} u_2 m_{2}^{-1}) \in K, \; \; 
\mbox{and} \; \; (h_1 u_1 h_{1}^{-1}, \theta_\alpha
(h_1 u_1 h_{1}^{-1})) \in K.
\]
Thus $m_{2} u_2 m_{2}^{-1} \theta_\alpha (h_1 u_1 h_{1}^{-1})^{-1}
\in X_2 \subset Z_{A_2}$, so $m_{2} u_2 m_{2}^{-1} \in H_2U_{2}^{\beta}$.
Thus $m_2$ normalizes both Borel 
subgroups of $M_{A_2}$ defined by $\De^+_{A_2}$
and by $-\De_{A_2}^+$. This implies $m_2 \in H_2$ which shows
that $K(H_1) \subset H_2$. Similarly, one shows that
$K(H_2) \subset H_1$.

Let now $D_1$ be any subset of $A_1$.
Assume that $(m_1, m_2) \in (M_{D_1} 
\times M_{A_2}) \cap K$. Write $m_1 = h_1 m_{1}^{\prime}$,
where $h_1 \in H_1$ and $m_{1}^{\prime} \in M_{D_1}^{\prime}$.
Since $M_{D_1}^{\prime}$ is generated by the $U_{1}^{\alpha}$'s
for $\alpha \in \De_{D_1}$, there exists 
$m_{2}^{\prime} \in M_{a(D_1)}^{\prime}$
such that $(m_{1}^{\prime}, m_{2}^{\prime}) \in K$. Thus
$(h_1, m_2 (m_{2}^{\prime})^{-1}) \in K$. It follows from
$K(H_1) \subset H_2$ that $m_2 (m_{2}^{\prime})^{-1} \in H_2$.
Thus $m_2 \in H_2 M_{a(D_1)}^{\prime} = M_{a(D_1)}$.
This shows that $K(M_{D_1}) \subset M_{a(D_1})$. One
proves similarly that $K(M_{a(D_1)}) \subset
M_{D_1}$. 
It  is clear from the definition that
$(M_{D_1} \times M_{a(D_1)}) \cap K$ is a generalized 
$a|_{D_1}$-graph. This proves 1).

Let $\theta \colon K_1/X_1 \to K_2/X_2$
be given as in \eqref{theta}. Since $X_1 \cap U_{D_1}^{A_1} =
\{e\}$ and $X_2 \cap U_{a(D_1)}^{A_2} = \{e\}$, we can
embed $U_{D_1}^{A_1}$ into $K_1/X_1$ and 
$U_{a(D_1)}^{A_2}$ into $K_2/X_2$. Then $\theta$ induces
a group isomorphism from $U_{D_1}^{A_1}$ to
$U_{a(D_1)}^{A_2}$ whose graph is
the intersection $(U_{D_1}^{A_1} \times U_{a(D_1)}^{A_2}) \cap K$.
This proves 2). 
 
To prove \eqref{PK1} and \eqref{PK2}, assume that $(p_1, m_2) \in 
(P_{D_1} \times M_{A_2}) \cap K$. Write $p_1 = m_1 u_1$,
where $m_1 \in M_{D_1}$ and $u_1 \in U_{D_1}$. Note that
$u_1 \in 
U_{D_1}^{A_1}$ because
$p_1, m_1 \in M_{A_1}$.
By 2), there exists $u_2 \in U_{a(D_1)}^{A_2}$ such that
$(u_1, u_2) \in K$. Thus $(m_1, m_2 u_{2}^{-1}) \in K$.
By 1), $m_2 u_{2}^{-1} \in M_{a(D_1)}$. Thus
$m_2 \in M_{a(D_1)}U_{a(D_1)}^{A_2} \in P_{a(D_1)}$, and
\[
(p_1, m_2) = (m_1, m_2u_{2}^{-1})(u_1, u_2) \in 
\left((M_{D_1} \times M_{a(D_1)})\cap K\right) 
\left((U_{D_1}^{A_1} \times U_{a(D_1)}^{A_2}) \cap K \right).
\]
This shows that 
\[
(P_{D_1} \times M_{A_2}) \cap K = 
(P_{D_1} \times P_{a(D_1)} ) \cap K =
\left((M_{D_1} \times M_{a(D_1)})
\cap K\right) \left((U_{D_1}^{A_1} \times U_{a(D_1)}^{A_2}) 
\cap K \right).
\]
Similarly, one proves that
$(M_{A_1} \times P_{a(D_1)}) \cap K$ is also equal to any
one of the above three groups. This proves \eqref{PK1} and \eqref{PK2}.

The identities in \eqref{PR1} follow directly from those in \eqref{PK1}.
To prove \eqref{PR2} and \eqref{PR3},
assume now that $(p_1, p_2) \in (P_{D_1} \times P_{a(D_1)}) \cap R_\A$.
Write $p_1 = m_1u_1$ and $p_2 = m_2 u_2$, where
$m_1 \in M_{D_1}, u_1 \in U_{D_1}, m_2 \in M_{a(D_1)}$ and
$u_2 \in U_{a(D_1)}$. Further write
\[
u_1 = u_{D_1}^{A_1}u_{A_1}, \; \; \;
u_2 = u_{a(D_1)}^{A_2}u_{A_2}, 
\] 
where
\[
u_{D_1}^{A_1} \in U_{D_1}^{A_1}, 
\; \; \;u_{A_1} \in U_{A_1}, \; \; \;
u_{a(D_1)}^{A_2}  \in u_{a(D_1)}^{A_2}, \; \; \;u_{A_2}
\in U_{A_2}.
\]
Since $(u_{A_1}, u_{A_2}) \in R_\A$, we have
\[
\left(m_1, u_{D_1}^{A_1}, \, 
m_2 u_{a(D_1)}^{A_2}\right) \in  \left(P_{D_1}^{A_1} 
\times P_{a(D_1)}^{A_2} \right) \cap R_\A = \left(P_{D_1}^{A_1} 
\times P_{a(D_1)}^{A_2} \right) \cap K.
\]
By \eqref{PK2}, we see that $(m_1, m_2) \in R_\A$. This shows 
\eqref{PR2}. Note that $(M_{D_1} \times M_{a(D_1)})\cap R_\A
=(M_{D_1} \times M_{a(D_1)})\cap K$. It is also easy to show that
\[
(U_{D_1} \times U_{a(D_1)}) \cap R_\A =  
\left((U_{D_1}^{A_1} \times U_{a(D_1)}^{A_2}) \cap K \right)
\left(U_{A_1} \times U_{A_2} \right).
\]
Now \eqref{PR3} follows from \eqref{PK2}.
\end{proof}

Finally we use \leref{quintuple} to treat compositions
of generalized graphs.

\ble{composition} Let $G_1,$ $G_2$ and $G_3$ be 
connected reductive 
algebraic groups (with fixed choices of maximal tori and 
Borel subgroups). Assume that $(a, K)$ and $(c, L)$
are two admissible pairs for $G_1 \times G_2$ and 
$G_2 \times G_3,$ respectively, such that the domain of $c$ 
coincides with the range of $a.$ Then the composition 
$(ca, L \circ K)$ is an admissible pair for $G_1 \times G_3.$
\ele

\begin{proof} \leref{quintuple} implies that 
for any root $\alpha$ of $G_1$ in the span of the 
domain of $a$, there exist group isomorphisms 
$\theta \colon U_1^{\alpha} \ra U_2^{a(\alpha)}$ and
$\varphi \colon U_2^{a(\alpha)} \ra U_3^{ca(\alpha)}$ 
such that
\begin{align*}
K \cap (U_1^{\alpha} \times G_2) &=
(\{e\} \times X_2) \,
(\Id \times \theta)(U_1^{\alpha}) 
\quad
\mbox{and} \\
L \cap (G_2 \times U_3^{ca(\alpha)}) &=
(Y_2 \times \{e\} ) \,
(\Id \times \varphi)(U_3^{a(\alpha)})
\end{align*}
for some subgroups $X_2$ and $Y_2$ of the fixed 
maximal torus of $G_2.$ Then
\[
(L \circ K) \cap 
(U_1^{\alpha} \times U_3^{ca(\alpha)}) =
(\Id \times \varphi \theta) (U_1^{\alpha}).
\] 
Therefore $L \circ K$ is a generalized $ca$-graph and
$(ca, L \circ K)$ is an admissible pair for $G_1 \times G_3.$
Let us note that it is not hard to determine the 
$ca$-quintuple corresponding to the composition
$L \circ K,$ cf. \leref{quintuple}, 
but it will not be needed and will omit it.
\end{proof}

\sectionnew{Some facts on Weyl groups and intersections of
parabolic subgroups}
\lb{facts}

\subsection{The Bruhat Lemma}
\lb{sub-Bruhat}
Fix a connected reductive algebraic group $G$ over $\k$ with a maximal 
torus $H$ and a set $\Delta^+$ of positive  roots for $(G, H)$.
Let $\Ga$ be the set of simple roots in $\Delta^+$. We will denote 
by $W$ the Weyl group of $(G, H)$ and by $l(\cdot)$ the standard
length function on $W.$ For $w \in W,$ $\dot{w}$ will denote a 
representative of $w$ in the normalizer of $H$ in $G$. 

Given $A \subset \Ga$, we will denote by $W_A$ the subgroup of
the Weyl group $W$ generated by elements in $A$. For $A \subset \Ga$,
let $P_A$ be the corresponding parabolic subgroup of $G$ containing 
the Borel subgroup of $G$ determined by $\Delta^+$, and let $M_A$ and $U_A$ 
be respectively its Levi factor containing $H$
and its unipotent radical.   For $A, C 
\subset \Ga$, the standard
Bruhat decomposition for $(P_A, P_C)$ double cosets of $G$ states
that 
\begin{equation}
\lb{eq_G-Bruhat}
G = \coprod_{[w] \in W_A \backslash W /W_C} P_A
\dot{w} P_C.
\end{equation}
We note that if $\wt{P}_A$ is a subgroup of $P_A$ such 
that $P_A = \wt{P}_A H$, then we also have
\begin{equation}
\lb{eq_G-Bruhat-Pa}
G = \coprod_{[w] \in W_A \backslash W /W_C} \wt{P}_A
\dot{w} P_C,
\end{equation}
which is a fact that will be used later in the paper.
\subsection{Minimal length representatives of double cosets}
\lb{sub-minimal}
Let $A, C \subset \Ga$.
Each $(W_A, W_C)$ double coset of $W$ contains a unique
element of minimal length, see e.g. \cite[Proposition 2.7.3]{Car}.
The set of minimal length representatives for 
$(W_A, W_C)$ double cosets will be denoted by
${}^A \! W^C.$ The latter set consists of
exactly those elements $w \in W$ with the 
properties 
\begin{equation}
\lb{eqw}
w^{-1}(A) \subset \De^+ \; \mbox{and} \;
w(C) \subset \De^+,
\end{equation}
see e.g. \cite[$\S$2.7]{Car}.
 When $A$ is empty, we set ${}^{A} \! W^C = W^C$.
If $D$ and $E$ are two subsets of $A \subset \Ga$, the set of
minimal length representatives 
in $W_A$ for the double cosets from 
$W_D \backslash W_A /W_E$ will be denoted by 
${}^D \! W_{A}^E$. Elements in ${}^D \! W_{A}^E$  
will be thought of as 
elements of $W$ via the inclusion of $W_A$ in $W.$
 
\bpr{W_prod} Fix three subsets $A$, $D$,  and $C$ of $\Ga$ with
$D \subset A$. Then 
every element $v \in {}^D\!W^C$ is represented in a unique way
as a product 
\begin{equation}
\lb{w_prod1}
v = u w \; \; \mbox{for some}
\; \;
w \in {}^{A} \! W^C, \;
u \in {}^D \!W_A^{A \cap w (C)}, 
\end{equation}
and all such products belong to 
${}^D\!W^C.$ Moreover, $l(v) = l(u) + l(w)$ for $u, v$, and $w$ in 
\eqref{w_prod1}.
\epr
Note that the defining set for $u$ in \eqref{w_prod1} depends on the 
first element $w.$ A proof of this result can be found in 
\cite[Proposition 2.7.5]{Car}, see also \cite[Lemma 4.3]{Y}.

By taking inverses in \eqref{w_prod1}, one sees that for three subsets
$A$, $E$, and $C$ of $\Ga$ with $E \subset C$,  every element 
$v \in {}^A \! W^E$ is represented in a unique way
as a product  
\begin{equation}
\lb{w_prod2}
v = w u \; \mbox{for some}
\; 
w \in {}^A \! W^C, \;
u \in {}^{w^{-1}(A) \cap C} \! W_C^E, 
\end{equation}
and all such products belong to ${}^A \! W^D.$
We also note that if $w \in {}^A \! W^C$, the minimality conditions
on $w$ imply that  
\begin{equation}
\lb{cap}
\De_A^+ \cap w( \De_C^+) = \De_{A \cap w(C)}^+ \hspace{.2in} 
\mbox{and} \hspace{.2in} A \cap w( \De_C^+) = A \cap w(C).
\end{equation}
As in $\S$\ref{adm},  we set 
$\De_D^+ = \Delta^+ \cap {\mathbb Z}[D]$ for $D \subset \Ga$.
\subsection{Intersections of parabolic subgroups}
\lb{sub-int-par} For $A \subset \Gamma$ and $w \in W$, we set 
\begin{equation}
\lb{adw}
w(P_A) = \Ad_{\dot{w}}(P_A), \; w(M_A) = \Ad_{\dot{w}}(M_A), \;
w(U_A) = \Ad_{\dot{w}}(U_A)
\end{equation}
all of which are independent of the choice of the representative $\dot{w}$. 
We will use $\pi_A$ to denote the
projection $P_A \to M_A$ with respect to the decomposition
$P_A = M_A U_A$.
For a subset $D$ of $A$, we will set $P_{D}^{A} = P_D \cap M_A$
and $U_{D}^{A} = U_D \cap M_A$.

\ble{int} Assume that $A, C \sub \Ga$ and $w \in {}^A \! W^C$. 
Then 

1) we have the direct product decompositions 
\begin{eqnarray}
\lb{PP1}
P_A \cap w(P_C) &=& (M_A \cap w(M_C)) \,(M_A \cap w(U_C))\, 
(U_A \cap w(M_C))\,
(U_A \cap w(U_C)) \\
\lb{PP2}
& = &M_{A \cap w(C)} \, U^A_{A \cap w(C)} \,
\left(w( U^C_{w^{-1}(A) \cap C})\right) \, (U_A \cap w(U_C))\\
\lb{PP3}
&=&M_{A \cap w(C)} \left(U_{A \cap w(C)} \cap w(U)\right)\\
\lb{PP4}
&= &M_{A \cap w(C)} \left(U_{A \cap w(C)} \cap w(U_{w^{-1}(A)
\cap C})\right)
\end{eqnarray} 
where $M_A \cap w(M_C)=M_{A \cap w(C)}$ normalizes all terms.
In particular, 
$P_A \cap w(P_C) \sub P_{A \cap w(C)}$ and 
$\pi_A(P_A \cap w(P_C)) = P_{A \cap w(C)}^{A};$

2) the following dimension formulas hold for the groups in 1):
\begin{align*}
&\dim U^A_{A \cap w(C)} + \dim (U_A \cap w(U_C))
= \dim U_C - l(w),
\\
& \dim w( U^C_{w^{-1}(A) \cap C}) +
\dim (U_A \cap w(U_C))
= \dim U_A - l(w).
\end{align*}
\ele

\begin{proof} The first part of \leref{int} is 
Theorem 2.8.7 from \cite{Car}. 
It only uses \eqref{eqw} and \eqref{cap}. 

For Part 2), we only prove the second formula since the proof of 
the first one is similar. To this end, note that  
$w( U^C_{w^{-1}(A) \cap C}) (U_A \cap w(U_C))\subset U_A$ 
and that it is the product of the 
one parameter unipotent subgroups of $G,$ corresponding to roots
\[
\al \in \De^+ - \De_A^+ \; \mbox{such that} \; w^{-1}(\al) \in \De^+.
\]
Since $w \in {}^A \! W$,  one has 
$w^{-1}(\al) \in \De^+$ for any $\al \in \De_A^+$.  
Thus the codimension in $U_A$ of 
$w( U^C_{w^{-1}(A) \cap C}) (U_A \cap w(U_C))$  is equal to the 
cardinality of the set $\De^+ \cap w(-\De^+)$  which is equal to 
the length of $w,$ and hence the second formula in 2). 
\end{proof}
\subsection{A lemma on orbits of subgroups}
The following lemma, which will be used later in the paper in conjunction 
with the Bruhat lemma, describes the orbit space of a group on
a set in terms of that of a bigger group.  We omit the proof
since it is straightforward. 

\ble{abs_ind}
Let $P$ be a group acting (on the left) on a set $M$ and let
$R$ be a subgroup of $P.$
Suppose that $\SS$ is a subset of $M$ parametrizing the
$P$ orbits on $M.$ For each $m \in \SS$, let $\Stab_m$ be the stabilizer
subgroup of $P$ at $m$. Then every $R$ orbit in $M$ contains
a point of the form $p\cdot m$ for a unique $m \in \SS$ and some $p \in P$.
Two points $p_1\cdot m$ and $p_2 \cdot m$, where $m  \in \SS$ and
$p_1, p_2 \in P$, are in the same $R$ orbit if and only if
$p_1$ and $p_2$ are in the
same $(R, \Stab_m)$ double coset in $P$.
\ele

\sectionnew{Proof of \thref{main1}}
\lb{prmain1}

\subsection{A special case}\lb{finstep}

For two arbitrary admissible quadruples $\A$ and $\C$, we
 will obtain  a description of 
$(R_\A, R_\C)$ double cosets in 
$G_1 \times G_2$ by an induction procedure to be presented
in $\S$\ref{ind-step}. In this section, 
we look at the last step in the induction  as
a special case. 

Assume that 
both $a$ and $c$ have domains $\Ga_1$ and ranges $\Ga_2$.  
Then $R_\A = K$ and $R_\C = L$. Modulo problems with the
centers, we can think of  both $K$ and $L$ as
graphs of isomorphisms from $G_1$ to $G_2$.
 
Denote the centers and the derived subgroups of $G_1$ and $G_2$ 
respectively by $Z_1, Z_2$, $G'_1$ and $G'_2$. Recall 
that $\eta_2 \colon G_1 \times G_2 \to G_2$ is
the projection to the second factor and that
$K_2 = \eta_2(K)$,  $L_2 = \eta_2(L)$. Let
$Z(K_2)$ and $Z(L_2)$ be the centers of $K_2$ and $L_2$ 
respectively. It follows from $K_2, L_2 \supset G_{2}^{\prime}$ 
that $Z(K_2) = Z_2 \cap K_2$ and $Z(L_2) = Z_2 \cap L_2$. Set 
\[
\Z_2 = Z(K_2)Z(L_2) \subset Z_2.
\]
For each $s \in Z_2/\Z_2$, fix a representative $\dot{s}$ of $s$ 
in $Z_2$. We will also introduce the group
\begin{equation}
\lb{J}
J =\{(k_1, l_1) \in G_1 \times G_1 \mid 
\exists \; g_2 \in G_2 \;
\mbox{such that} \; (k_1, g_2) \in K, (l_1, g_2) \in L\}.
\end{equation}
In terms of compositions of group correspondences in 
$\S$\ref{notation}, 
$J = \sigma(L)  \circ K$, where
$\sigma \colon G_1 \times G_2 \to G_2 \times G_1 \colon (g_1, g_2)
\mapsto (g_2, g_1)$.
By \leref{composition},   
$J$ is a generalized $c^{-1}a$-graph.
Let  $J$ act on $G_1$  from the left by
\begin{equation}
\lb{jact}
(k_1, l_1) \cdot g_1 = k_1 g_1 l_{1}^{-1}, \hspace{.2in} 
(k_1, l_1) \in J, \; g_1 \in G_1.
\end{equation}

\ble{end}
When both $a$ and $c$ have domains $\Ga_1$ and ranges $\Ga_2$, 
every  $(R_\A, R_\C)$ double coset of $G_1 \times G_2$ is of the form
\[
[g_1, \dot{s}] \; \mbox{for some} \; g_1 \in G_1, \;  s \in Z_2 /\Z_2.
\]
Two such cosets $[g_1, \dot{s}]$ and $[h_1, \dot{t}]$ are the 
same  
if and only if $s = t$ and $g_1$ and $h_1$ are in the
same orbit for the $J$ action on $G_1$ given by 
\eqref{jact}, cf. also \eqref{J}.
\ele

\begin{proof}  Using the fact that $G_2 = Z_2L_2$, we see that
every $(R_\A, R_\C)$ double coset of $G_1 \times G_2$ is of the form
$[g_1, z]$ for some $g_1 \in G_1$ and $z \in Z_2$. By writing
$z = \dot{s} z_K z_L$ for some  $s \in Z_2 /\Z_2$, 
$z_K \in Z(K_2)$, and $z_L \in Z(L_2)$,
and by changing $g_1$ to another element in $G_1$ if necessary, 
we see that
every $(R_\A, R_\C)$ double coset of $G_1 \times G_2$ is of the form
$[g_1, \dot{s}]$ for some $g_1 \in G_1$ and $s \in Z_2/\Z_2$.
It is clear that if $g_1, h_1 \in G_1$ are in the same $J$ orbit
in $G_1$, then $[g_1, \dot{s}] = [h_1, \dot{s}]$ for any 
$s \in Z_2/\Z_2$. Assume now that $[g_1, \dot{s}]=[h_1, \dot{t}]$ for
some $g_1, h_1 \in G_1$ and $s, t \in Z_2/\Z_2$. Then there exist
$(k_1, k_2) \in K, (l_1, l_2) \in L$
such that 
\[
(h_1, \dot{t})= (k_1, k_2) (g_1, \dot{s}) (l_1, l_2)^{-1}=
(k_1 g_1 l_{1}^{-1}, \, \dot{s} k_2 l_{2}^{-1}).
\]
We claim that
$k_2 l_{2}^{-1} \in \Z_2$. Indeed, write $k_2 = k'_2 z_K$,
where $k'_2 \in G'_2$ and $z_K \in Z_2$. Since $G'_2 \subset
K_2$, $z_K \in K_2 \cap Z_2 = Z(K_2)$.
Similarly, we can write $l_{2}^{-1}= l'_2 z_L$, where
$l'_2 \in G'_2$ and $z_L \in Z(L_2)$. Thus
$k_2 l_{2}^{-1} = (k'_2 l'_2) z_K z_L$.
On the other hand, $k_2 l_{2}^{-1} = \dot{t} \dot{s}^{-1} \in Z_2$,
so $k'_2 l'_2 \in G'_2 \cap Z_2 \subset Z(K_2)$.
Thus $k_2 l_{2}^{-1} = (k'_2 l'_2 z_K) z_L
\in \Z_2$. It now follows that $\dot{t} \in \dot{s} \Z_2$,
so $s = t$, and $\dot{s} = \dot{t}$. Thus 
$k_2 = l_2$, and $(k_1, l_1) \in J$. Hence $g_1$ and $h_1$
are in the same $J$ orbit in $G_1$.
\end{proof}

\subsection{The main induction step in the proof of \thref{main1}}
\lb{ind-step}

In this section, we will reduce the classification of $(R_\A, R_\C)$ 
double cosets in $G_1 \times G_2$ to that of similar double cosets
in $M_{A_1} \times M_{C_2}$.   
More precisely, for each pair $(w_1, w_2) \in {}^{A_1}\!W_1^{C_1} 
\times {}^{A_2}\!W_2^{C_2}$, 
we will define two
new admissible pairs 
\[
\A^\newnew(w_1, w_2) \; \; \; \mbox{and}
\; \; \; 
\C^\newnew(w_1, w_2)
\]
for $M_{A_1} \times M_{C_2}$.
Denoting by $\RAnew$  and 
$\RCnew$ the subgroups 
of $M_{A_1} \times M_{C_2}$
defined according to \eqref{R}, corresponding respectively to the 
admissible pairs $\A^\newnew(w_1, w_2)$ and $\C^\newnew(w_1, w_2)$, we will
show that  
the double coset space $R_\A \backslash G_1 \times G_2/R_\C$
can be identified with the union
over ${}^{A_1}\!W_1^{C_1} \times {}^{A_2}\!W_2^{C_2}$ of 
the double coset spaces
\[
\RAnew \backslash M_{A_1} \times M_{C_2}/\RCnew.
\]
 
Let   
$(w_1, w_2) \in {}^{A_1}\!W_1^{C_1} \times {}^{A_2}\!W_2^{C_2}$.
Set
\begin{eqnarray}
\lb{Anew}
\Aanew &= a^{-1}(A_2 \cap w_2(C_2)), \hs 
\Abnew = C_2 \cap w_{2}^{-1}(A_2) \\
\lb{Cnew}
\Canew &= A_1 \cap w_1(C_1), \hs 
\Cbnew = c(C_1 \cap w_{1}^{-1}(A_1)).
\end{eqnarray}
We will regard
\[
w_{2}^{-1}a \colon \Aanew \rightarrow \Abnew
\; \; {\rm and}  \; \; \;  \\
c w_{1}^{-1} \colon \Canew \rightarrow \Cbnew
\]
as partial isometries from $A_1$ to $C_2$. 
For $i = 1, 2$, let $\dot{w}_i$ be a representative of $w_i$ in the 
normalizer of $H_i$ in $G_i$. Define
\begin{eqnarray}
\lb{Knew}
K^\newnew \new&= \left(M_{\Aanew} \times 
M_{\Abnew}\right) \cap \Ad_{(e, \dot{w}_2)}^{-1}K\\
\lb{Lnew}
L^\newnew \new &= \left(M_{\Canew} \times M_{\Cbnew}\right)
\cap \Ad_{(\dot{w}_1, e)} L.
\end{eqnarray} 
Then it follows from \leref{property} that
\[
\A^\newnew \new := (w_{2}^{-1}a, \; K^\newnew \new) \; \; \; 
\mbox{and} \; \; \; 
\C^\newnew \new := (c w_{1}^{-1}, \; L^\newnew \new)
\]
are admissible pairs for $M_{A_1} \times M_{C_2}$. Let
\begin{align*} 
\RAnew &= K^\newnew \new  \left(U_{\Aanew}^{A_1} \times
U_{\Abnew}^{C_2} \right) \\
\RCnew &=  
L^\newnew\new \left(U_{\Canew}^{A_1} \times U_{\Cbnew}^{C_2}
\right).
\end{align*}
  
\bpr{IndSt} {\bf{(Induction Step):}} 
Every $(R_{\A}, R_{\C})$ double
coset in $G_1 \times G_2$ is of the form
\[
[m_{A_1} \dot{w}_1, \dot{w}_2 m_{C_2}] \; \; \mbox{for some}
\; \;  w_i \in {}^{A_i} \! W_i^{C_i}, \;   
m_{A_1} \in M_{A_1}, \; m_{C_2} \in M_{C_2}.
\]
Two such cosets $[m_{A_1} \dot{w}_1, \dot{w}_2 m_{C_2}]$ and
$[m'_{A_1} \dot{w}'_1, \dot{w}'_2 m'_{C_2}]$  coincide if and only 
if $w'_1 = w_1, w'_2 = w_2,$ and  $(m_{A_1}, m_{C_2})$ and
$(m'_{A_1}, m'_{C_2})$ are in the same
$(\RAnew, \RCnew)$ double
coset of $M_{A_1} \times M_{C_2}.$ 
\epr

\begin{proof}
Recall that for $i = 1, 2$, 
\[
K_i = \eta_i (K) \subset M_{A_1} \hs \mbox{and} \hs
L_i = \eta_i(L) \subset M_{C_i}, 
\]
where $\eta_i \colon G_1 \times G_2 \to G_i$ is the projection to the
$i$'th factor. Let
\begin{align*}
\tP_{A_2} &= K_2U_{A_2} \subset P_{A_2}, \hs  
P_{\A, 1} = \left( M_{A_1} \times \{e\} \right) R_{\A} 
=P_{A_1} \times \tP_{A_2}\\
\tP_{C_1} &= L_1 U_{C_1} \subset P_{C_1}, \hs 
P_{\C, 2} = \left( \{ e \} \times M_{C_2} \right) R_{\C} =
\tP_{C_1} \times P_{C_2}.
\end{align*}
Consider the (left) action of $P_{\A, 1} \times P_{\C, 2}$
on $G_1 \times G_2,$ defined by
\begin{equation}
\lb{rac}
(p_{A_1},\tp_{A_2},\tp_{C_1},p_{C_2}) \cdot
(g_1,g_2)=
(p_{A_1} g_1 \tp_{C_1}^{-1}, \tp_{A_2} g_2 p_{C_2}^{-1})
\end{equation}
for $(p_{A_1},\tp_{A_2},\tp_{C_1},p_{C_2}) \in P_{\A, 1}
\times P_{\C, 2} $ and $(g_1, g_2) \in G_1 \times G_2.$
The set of orbits for this action coincides with
the set of $(P_{\A, 1}, P_{\C, 2})$ double cosets of $G_1 \times G_2$,
which, by \eqref{eq_G-Bruhat-Pa}, consists of the double cosets of the form
$P_{\A,1} (\dot{w}_1, \dot{w}_2) P_{\C,2}$
for 
$(w_1, w_2) \in  {}^{A_1} \! W_1^{C_1} \times {}^{A_2} \! W_2^{C_2}$.
We will apply \leref{abs_ind} for $P=P_{\A, 1} \times P_{\C, 1}$
and $R=R_{\A} \times R_{\C}$ to classify
$(R_{\A}, R_{\C})$ double cosets in $G_1 \times G_2$.
For $w_i \in {}^{A_i} \! W_i^{C_i}$ with $i = 1, 2$, let 
$\Stab_{(\dot{w}_1, \dot{w}_2)} \subset P_{\A, 1} \times P_{\C, 2}$
be the stabilizer subgroup of $P_{\A, 1} \times P_{\C, 2}$
at $(\dot{w}_1, \dot{w}_2)$.    In view of \leref{abs_ind}, to classify
$(R_{\A}, R_{\C})$ double cosets in $G_1 \times G_2$, it suffices to 
understand the space of double cosets
\begin{equation}
\lb{SRdouble}
(R_{\A} \times R_{\C})
\backslash
P_{\A, 1} \times P_{\C, 2}
/ \Stab_{(\dot{w}_1, \dot{w}_2)}
\end{equation}
for each pair  $(w_1, w_2) \in \WW$. To this end, consider the projection
\[
\pi:=\pi_{A_1} \times \pi_{A_2} \times \pi_{C_1} \times \pi_{C_2} \colon
P_{\A, 1} \times P_{C, 2} \rightarrow 
M_{A_1} \times K_2 \times L_1 \times M_{C_2}.
\]
(See $\S$\ref{facts} for notation). 
Since $M_{A_1} \times K_2 \times L_1 \times M_{C_2}$
normalizes $U_{A_1} \times U_{A_2} \times U_{C_1} \times U_{C_2}$,
$\pi$ induces an identification 
\[
(R_{\A} \times R_{\C})
\backslash
P_{\A, 1} \times P_{\C, 2}
/ \Stab_{(\dot{w}_1, \dot{w}_2)} \cong 
(K \times L) \backslash M_{A_1} \times K_2 \times L_1 \times M_{C_2}
/\pi(\Stab_{(\dot{w}_1, \dot{w}_2)}).
\]
Note now that we can identify 
\begin{equation}
\lb{QtoX}
(K\times L) \backslash (M_{A_1} \times K_2 \times L_1 \times M_{C_2})
\cong (X_1 \times Y_2)\backslash (M_{A_1} \times M_{C_2}),
\end{equation}
where $X_1 = \eta_1(\ker(\eta_2|_K)) \subset Z_{A_1}$ and 
$Y_2 = \eta_2(\ker(\eta_1|_L)) \subset Z_{C_2}$.
Indeed, for 
\[
(m_{A_1}, k_2, l_1, m_{C_2}) \in 
M_{A_1} \times K_2 \times L_1 \times M_{C_2},
\]
let 
$k_1 \in K_1$ and $l_{2} \in L_2$ be such that
$(k_1, k_2) \in K$ and $(l_1, l_{2}) \in L$.
Then
\[
(K \times L)(m_{A_1}, k_2, l_1, m_{C_2}) =
(K \times L)(k_1^{-1} m_{A_1}, \, e, \, e, \,  
l_{2}^{-1} m_{C_2}).
\]
It is easy to see that the map
\[
(K \times L)(m_{A_1}, k_2, l_1, m_{C_2}) \mapsto
(X_1 \times Y_2) (k_1^{-1} m_{A_1},  \,  
l_2^{-1} m_{C_2})
\]
gives a well-defined identification of the two spaces in 
\eqref{QtoX}, and that the right translation of
$(m_{A_1}, k_2, l_1, m_{C_2})$  on 
$(K \times L)\backslash 
(M_{A_1} \times K_2 \times L_1 \times M_{C_2})$
becomes the following map on $(X_1 \times Y_2)\backslash 
(M_{A_1} \times M_{C_2})$:
\begin{equation}
\lb{4M2M}
(X_1 \times Y_2) (m_1, m_2) 
\stackrel{(m_{A_1}, k_2, l_1, m_{C_2})}{\longmapsto}
(X_1 \times Y_2) (k_1^{-1}m_1 m_{A_1},
\, \, l_{2}^{-1}m_2 m_{C_2}).
\end{equation} 
Thus, every $(R_\A, R_\C)$ double coset in $G_1 \times G_2$
is of the form $[m_{A_1}\dot{w}_1, \dot{w}_2m_{C_2}^{-1}]$ 
for a unique pair
$(w_1, w_2) \in \WW$
and for some $(m_{A_1}, m_{C_2}) \in M_{A_1} \times M_{C_2}$. Two such
double cosets $[m_{A_1}\dot{w}_1, \dot{w}_2m_{C_2}^{-1}]$ and 
$[l_{A_1}\dot{w}_1, \dot{w}_2l_{C_2}^{-1}]$ coincide if 
and only if
$(X_1 \times Y_2)(m_{A_1}, m_{C_2})$ and 
$(X_1 \times Y_2)(l_{A_1}, n\l_{C_2})$
are in the same $\pi(\Stab_{(\dot{w}_1, \dot{w}_2)})$ orbit 
for the action of
$\pi(\Stab_{(\dot{w}_1, \dot{w}_2)})$ on 
$(X_1 \times Y_2)\backslash (M_{A_1} \times M_{C_2})$ given in
\eqref{4M2M}. 
Consider the map
\[
\nu \colon M_{A_1} \times M_{C_2} \ra
(X_1 \times Y_2)\backslash (M_{A_1} \times M_{C_2}) \colon
\nu(m_{A_1}, m_{C_2}) = 
(X_1 \times Y_2)( m_{A_1},  m_{C_2}^{-1}).
\]
The proof of \prref{IndSt} will be complete once we 
prove the following
fact about $\nu$: for every $(m_{A_1} \times m_{C_2}) 
\in M_{A_1} \times M_{C_2}$,
\begin{equation}
\lb{nu}
\nu^{-1} \left(
(X_1  \times Y_2)(m_{A_1},  m_{C_2}^{-1}) \cdot
\pi(\Stab_{(\dot{w}_1, \dot{w}_2)} )\right) =
\RAnew (m_{A_1}, m_{C_2}) \RCnew, 
\end{equation}
where again the action of $\pi(\Stab_{(\dot{w}_1, \dot{w}_2)})$ on 
$(X_1 \times Y_2)\backslash (M_{A_1} \times M_{C_2})$ is given in
\eqref{4M2M}. To prove this property of $\nu$, we set
\[
D_1 = A_1 \cap w_1(C_1), \; D_2 = A_2 \cap w_2(C_2), \; 
E_1 = C_1 \cap w_{1}^{-1}(A_1),\; E_2 = C_2 \cap w_{2}^{-1}(A_2).
\]
Then we know from 1) of \leref{int} that 
$\pi(\Stab_{(\dot{w}_1, \dot{w}_2)})$
consists precisely of all the elements of the form
\[
(\Ad_{\dot{w}_1} (l_1)u_{D_1}^{A_1}, \,\; k_2
u_{D_2}^{A_2}, \,\; l_1 u_{E_1}^{C_1}, \, \;
\Ad_{\dot{w}_2}^{-1} (k_2) u_{E_2}^{C_2}),
\]
where $u_{D_i}^{A_i} \in U_{D_i}^{A_i}$, 
$u_{E_i}^{C_i} \in U_{E_i}^{C_i}$ for $i = 1, 2$, 
$l_1 \in L_1 \cap M_{E_1}$,
and
$k_2 \in K_2 \cap M_{D_2}$. 
It is then easy to see from the definition 
of the action in \eqref{4M2M} and from
the properties of the groups $K$ and $L$ as stated in 2) and 3) of
\leref{property} that 
\eqref{nu} holds.
\end{proof}

\subsection{Proof of \thref{main1}}\lb{proof}
  
Fix two elements $v_1 \in W_1^{C_1}$ and $v_2 \in {}^{A_2} \! W_2.$
According to \prref{W_prod} and \eqref{w_prod2}, they can
be uniquely decomposed as
\begin{align}
\lb{v1u1}
&v_1 = u_1 w_1,\; 
w_1 \in {}^{A_1} \! W_1^{C_1}, \;
u_1 \in W_{A_1}^{A_1 \cap w_1 (C_1)}; 
\\ 
\lb{v2u2}
&v_2 = w_2 u_2, \; 
w_2 \in {}^{A_2} \! W_2^{C_2}, \;
u_2 \in {}^{w_2^{-1}(A_2) \cap C_2} \! W_{C_2}. 
\end{align}
Recall the definition of   $A^{\newnew}_{1}\new$  in
\eqref{Anew}. 
Set
\begin{align}
\lb{Astabnew}
\Aanew(u_1, u_2) =\{ \al \in \Aanew \mid
&(u_1 (c w_1^{-1})^{-1} u_2^{-1} (w_2^{-1}a))^n \al \; 
\mbox{is well defined}
\\ 
\nn
&\mbox{and is in} \; \Aanew \; \mbox{for} \; n=0, 1, \ldots\}.
\end{align}
Recall that $A_1(v_1, v_2)$ is defined in \eqref{A1stab}.

\ble{ind_stab} Let $v_1 \in W_1^{C_1}$ and $v_2 \in {}^{A_2} \! W_2$
be decomposed as $v_1 = u_1 w_1$ and $v_2 = w_2 u_2$ as in 
\eqref{v1u1}-\eqref{v2u2}. Then $\Aanew(u_1, u_2) = A_1(v_1, v_2)$.
\ele

\begin{proof} Note that the composition of maps in \eqref{Astabnew}
is nothing but
\[
(u_1 w_1) c^{-1} (u_2^{-1} w_2^{-1}) a = v_1 c^{-1} v_2^{-1} a.
\]
Thus if $\al \in \Aanew(u_1, u_2)$ then 
$(v_1 f^{-1} v_2^{-1} d)^n\al$ is well-defined and is a root of
$\Aanew \sub A_1.$ Therefore $\al \in A_1(v_1, v_2).$
Next we prove the opposite inclusion. First we show that 
if $\be \in A_1$ and $c^{-1} v_2^{-1} a(\be)$ is well-defined, then
$\be \in \Aanew.$ The assumptions imply that 
$a(\be) \in A_2$ and $u_2^{-1} w_2^{-1} a(\be) \in C_2.$
Since $u_2 \in W_{C_2}$, we get that 
$w_2^{-1} a(\be) \in \De_{C_2}$. Because
$w_2 \in {}^{A_2} \! W_2^{C_2}$, we further deduce
that $w_2^{-1} a(\be) \in \De_{C_2}^+.$ Now
$a(\be) \in A_2 \cap w_2(\De_{C_2}^+)$ and due to
\eqref{cap}, $a(\be) \in A_2 \cap w_2(C_2),$
so $\be \in \Aanew.$
If $\al \in A_1(v_1, v_2)$ then 
$(c^{-1} v_2^{-1} a) ( v_1 c^{-1} v_2^{-1} a)^n (\al)$
needs to be well defined for all $n=0, 1, \ldots$ so
$\al \in \Aanew$ and $( v_1 c^{-1} v_2^{-1} a)^n (\al)$
is well defined and is in $\Aanew$ for all 
$n=1, 2, \ldots$ Therefore $\al \in \Aanew(u_1, u_2).$
\end{proof}

Recall the groups $\Kvv$ and $\Lvv$ given in \eqref{Kvv} 
and \eqref{Lvv}.
It is easy to see from  \leref{end} that \thref{main1} is 
equivalent to the following proposition.

\bpr{main2} Any $(R_\A, R_\C)$ double coset of $G_1 \times G_2$ 
is of the form $[m_1 \dot{v}_1, \dot{v}_2 m_2]$ for some 
$v_1 \in W_1^{C_1},$ $v_2 \in {}^{A_2} \! W_2$,  
$m_1 \in M_{A_1(v_1, v_2)},$ and
$m_2 \in M_{C_2(v_1, v_2)}$.  
Two such cosets
$[m_1 \dot{v}_1, \dot{v}_2 m_2]$ and
$[l_1 \dot{v}'_1, \dot{v}'_2 l_2]$ 
coincide if and only $v'_i = v_i$ for $i = 1, 2$,  and
$(m_1, m_2)$ and $(l_1, l_2)$  are in the same 
$(\Kvv, \Lvv)$ double coset in
$M_{A_1(v_1, v_2)} \times M_{C_2(v_1, v_2)}$.
\epr

\noindent
{\em{Proof of \prref{main2}.}}
It is not hard to see that by repeatedly using  
\prref{IndSt} and by using \leref{ind_stab}, 
we can reduce the description of 
$(R_\A, R_\C)$ double cosets in $G_1 \times G_2$ to that of
the special cases treated in \leref{end}.
\prref{main2} holds trivially for such cases. 
Thus it suffices to prove \prref{main2} for $G_1 \times G_2$ 
by assuming that it holds for $M_{A_1} \times M_{C_2}$. 

For
$(w_1, w_2) \in \WW$ 
we set 
\[
\bD(w_1, w_2)= \RAnew \backslash
M_{A_1} \times M_{C_2} /\RCnew
\]
and use $[m_{A_1}, m_{C_2}]'$ to denote the point in $\bD(w_1, w_2)$
defined by $(m_{A_1}, m_{C_2}) \in M_{A_1} \times M_{C_2}$.
Then by 
\prref{IndSt}, we have the disjoint union
\[
G_1 \times G_2 = \bigsqcup_{(w_1, w_2)} \; \; \; 
\bigsqcup_{[m_{A_1}, m_{C_2}]' \in \bD(w_1, w_2)} 
[m_{A_1} \dot{w}_1, \dot{w}_2 m_{C_2}],
\]
where $(w_1, w_2) \in {}^{A_1}\!W_1^{C_1} \times {}^{A_2}\!W_2^{C_2}$.
 For  such a pair $(w_1, w_2)$, and for 
$u_1 \in W_{A_1}^{A_1 \cap w_1(C_1)}$ and $u_2 \in 
{}^{C_2 \cap w_{2}^{-1}(A_2)}\!W_{C_2}$, set $v_1 = u_1w_1$ and
$v_2 = w_2u_2$. Then we know from \prref{W_prod} and \eqref{w_prod2}
that $v_1 \in W_1^{C_1}$ and $v_2 \in {}^{A_2}\!W_2$, and that
every $v_1 \in W_1^{C_1}$ and $v_2 \in {}^{A_2}\!W_2$ are of this form.
Choose representatives $\hat{u}_1$ and $\hat{u}_2$ of $u_1$ and $u_2$
in the normalizers of $H_1$ and $H_2$ in $M_{A_1}$ and $M_{C_2}$ 
respectively such that 
$\dot{v}_1 = \hat{u}_1\dot{w}_1$ and $\dot{v}_2 = \dot{w}_2\hat{u}_2$.
Set
\[
\bD(v_1, v_2) = \Kvv
\backslash M_{A_1(v_1, v_2)} \times M_{C_2(v_1, v_2)}/
\Lvv,
\]
and let
$[m_1, m_2]_{\mathrm{fin}}$ denote the point in $\bD(v_1, v_2)$ 
defined by
$(m_1, m_2) \in M_{A_1(v_1, v_2)} \times M_{C_2(v_1, v_2)}$.  
Then by applying
\prref{main2} to $\bD(w_1, w_2)$ and by using  \leref{ind_stab}, 
we know that 
\[
\bD(w_1, w_2) = \bigsqcup_{(u_1, u_2)} \; \; \; 
\bigsqcup_{[m_1, m_2]_{\mathrm{fin}}
\in \bD(v_1, v_2)} [m_1 \hat{u}_1,  \hat{u}_2 m_2]'
\]
is a disjoint union, where $(u_1, u_2) \in W_{A_1}^{A_1 \cap w_1(C_1)}
\times {}^{C_2 \cap w_{2}^{-1}(A_2)}\!W_{C_2}$. Thus we have 
\[
G_1 \times G_2 = \bigsqcup_{v_1 \in W_1^{C_1}, v_2 \in {}^{A_2}\!W_2}
\; \; \; \bigsqcup_{[m_1, m_2]_{\mathrm{fin}}
\in \bD(v_1, v_2)} [m_1  \dot{v}_1, \dot{v}_2 m_2]
\]
as a disjoint union.
\qed

\sectionnew{Proofs of \thref{main2} and \thref{dim}}
\lb{prmain2}

We will keep the notation as in $\S$\ref{classification}
and  $\S$\ref{stabilizer}. 

\subsection{Proof of \thref{main2}}\lb{prmain2-1}

\ble{stab}
For $q = (m_1 \dot{v}_1, \dot{v}_2 s_2)$, where 
$(v_1, v_2) \in W_{1}^{C_1} \times {}^{A_2}\!W_2$, $m_1  \in M_{\Aavv}$ 
and $s_2 \in Z_{\Cbvv}$, let
\begin{eqnarray}
\lb{MStab-a}
\MStab(q) &=& \Kvv \cap \Ad_{(m_1, e)} \Lvv \\
\lb{MStab-b}
& =&\left(M_{\Aavv} \times M_{\Cbvv} \right) \cap 
\Ad_{(e, \dot{v}_2)}^{-1}K \cap 
\Ad_{(m_1 \dot{v}_1, e)}L\\
\lb{UStab}
\UStab(q) & =& \Stab(q) \cap 
\left( (U_{\Aavv} \cap v_1(U_1)) \times 
(U_{\Cbvv} \cap v_2^{-1}(U_2))\right).
\end{eqnarray}
Then 
$\Stab(q) = \MStab(q) \; \UStab(q)$.
\ele

\begin{proof} 
By \prref{W_prod} and \eqref{w_prod2}, we can write
$v_1 = u_1 w_1$ and $v_2 = w_2 u_2$, where 
$w_1 \in {}^{A_1}\!W_{1}^{C_1}, \; 
u_1 \in W_{A_1}^{A_1 \cap w_1(C_1)}, \; 
w_2\in {}^{A_2}\!W_{1}^{C_2}$ and $u_2 \in 
{}^{C_2 \cap w_{2}^{-1}(A_2)}\!W_{C_2}$. 
We will first give a description of $\Stab(q)$ in terms of 
the groups $K^\newnew \new$ and $L^\newnew \new$ in 
\eqref{Knew} and \eqref{Lnew}. Set
$\dot{u}_1 = \dot{v}_{1} \dot{w}_{1}^{-1}  
\in M_{A_1}$ and $\dot{u}_2 = 
\dot{w}_{2}^{-1} \dot{v}_{2} \in M_{C_2}$. 

Let $(p_1, p_2) \in \Stab(q)$.  Then $\Ad_{\dot{v}_2s_2}(p_2) \in 
P_{A_2} \cap w_2(P_{C_2})$ because $\dot{u}_2 s_2 \in P_{C_2}$. By 
1) of \leref{int}, $\Ad_{\dot{v}_2s_2}(p_2) 
\in P_{A_2 \cap w_2(C_2)}$. It follows from
$(p_1, \Ad_{\dot{v}_2s_2}(p_2)) \in R_\A$ and \eqref{PR1} that
$p_1 \in P_{a^{-1}(A_2 \cap w_2(C_2))}$. Hence
\[
(p_1, \Ad_{\dot{v}_2s_2}(p_2)) \in 
R_\A \cap \left(P_{a^{-1}(A_2 \cap w_2(C_2))} \times 
(P_{A_2} \cap w_2(P_{C_2}) \right).
\]
By 1) of \leref{int} and 3) of \leref{property}, we see that
$(p_1, \Ad_{\dot{v}_2s_2}(p_2))$ belongs to 
the group 
\[
\left((M_{a^{-1}(A_2 \cap w_2(C_2))} \times
M_{A_2 \cap w_2(C_2)}) \cap K \right)
\left(U_{a^{-1}(A_2 \cap w_2(C_2))} \times 
w_2(U_{C_2 \cap w_{2}^{-1}(A_2)}) \right).
\]
Thus
\[
(p_1, \Ad_{\dot{u}_2s_2} (p_2) ) 
\in K^\newnew \new \left(U_{\Aanew} \times U_{\Abnew} \right),
\]
where $\Aanew$ and $\Abnew$ are given in \eqref{Anew}.
Similarly, with $\Canew$ and $\Cbnew$ given in \eqref{Cnew}, 
one shows that 
\[
(\Ad_{m_1\dot{u}_1}^{-1}(p_1), p_2 ) \in  
L^\newnew \new 
\left(U_{\Canew} \times U_{\Cbnew} \right).
\]
Using \leref{ind_stab} one proves inductively that
\begin{align*}
& (p_1, \Ad_{s_2}(p_2)) \in \Kvv 
\left(U_{\Aavv} \times U_{\Cbvv} \right)\\
&(\Ad_{m_1}^{-1}(p_1), p_2) \in \Lvv 
\left(U_{\Aavv} \times U_{\Cbvv} \right).
\end{align*}
Thus
\begin{equation}\lb{pp}
(p_1, p_2 ) 
\in \left(\Kvv \cap \Ad_{(m_1, e)} \Lvv \right)
\left(U_{A_1(v_1, v_2)} \times U_{C_2(v_1, v_2)} \right).
\end{equation} 
It is easy to see that 
\[
\MStab(q):=\Kvv \cap \Ad_{(m_1, e)} 
\Lvv \subset \Stab(q)
\]
and that \eqref{MStab-b} holds. 
Thus, if we set 
\[
\UStab(q) = \Stab(q) \cap 
\left(U_{A_1(v_1, v_2)} \times U_{C_2(v_1, v_2)} \right),
\] 
it remains to 
show that $\UStab(q)$ is also given as in \eqref{UStab}.
We note that if $(p_1, p_2) \in \UStab(q)$, then 
$\Ad_{m_1}^{-1} p_1 \in v_1(P_{C_1})$. 
On the other hand,
$\Ad_{m_1}^{-1}(p_1)\in P_{\Aavv}$ by \eqref{pp}. Thus $ 
\Ad_{m_1}^{-1}(p_1) \in 
P_{\Aavv} \cap v_1(P_{C_1})$, so $p_1 \in P_{\Aavv} \cap v_1(P_{C_1})$.
 Now 
$v_1 \in {}^{\Aavv}\!W_{1}^{C_1}$  
because $v_{1}^{-1}(\Aavv) \subset C_1 \subset
\Delta_{1}^{+}$, Thus by \eqref{PP3}
in \leref{int}, $p_1 \in  U_{\Aavv} \cap v_1(U_1)$. 
Similarly, $p_2 \in U_{\Cbvv} \cap v_{2}^{-1}(U_2)$. Thus $\UStab(q)$
is given by \eqref{UStab}.
\end{proof}

Recall that $J_{m_1}(v_1,v_2)$ denotes the stabilizer subgroup
of $\Jvv$ at $m_1 \in M_{A_1(v_1, v_2)}$ for the action given in
\eqref{Q1-M1}.
The group $\MStab(q)$ consists of all pairs 
$(k_1, k_2) \in \Kvv$ such that 
$(\Ad^{-1}_{m_1}k_1, k_2) \in \Lvv.$ 
Firstly this implies that 
the group 
$(\ker \eta_1|_L) \cap (\Id \times v_2^{-1})(\ker \eta_1|_K)$
lies in the center of $\MStab(q)=\Kvv \cap \Ad_{(m_1, e)}\Lvv.$ 
Secondly, for 
$(k_1, k_2) \in \MStab(q)$ we have that 
$(k_1, \Ad^{-1}_{m_1} k_1) \in J_{m_1}(v_1,v_2)$
because of \eqref{compJ}. Therefore 
\[
\mu \colon \MStab(q) \ra J_{m_1}(v_1,v_2), \quad 
\mu(k_1, k_2) = (k_1, \Ad^{-1}_{m_1} k_1)
\]
is a well defined homomorphism. The following Lemma 
describes explicitly $\MStab(q).$ 
Its proof is straightforward and we will omit it. 

\ble{ext} The following is a group exact sequence
\[
\{e\} \ra (\ker \eta_1|_L) \cap 
(\Id \times v_2^{-1}) (\ker \eta_1|_K) 
\hra \MStab(q) \ra J_{m_1}(v_1,v_2) \ra \{e\}.  
\]
\ele

The next Lemma will be used to describe the structure of 
$\UStab(q),$ stated in \thref{main2}.

\ble{idemp} The maps $\sigma_1$ and $\sigma_2$ given by 
\eqref{sigma1}-\eqref{sigma2} are well defined endomorphisms
of $U_{\Aavv} \cap v_1(U_1)$ and
$U_{\Cbvv} \cap v_2^{-1}(U_2),$
respectively. Moreover for some sufficiently large 
integer $k,$ $\sigma_i^{k+1} \equiv e$ for $i = 1, 2$.
\ele

\begin{proof} 
We will only prove the statements 
for $\sigma_1$ as the ones for $\sigma_2$ are
analogous. First we show that
\[
\til{\sigma}_1 = \Ad_{\dot{v}_1} \theta_{c}^{-1} \pi_{C_2}
\Ad_{\dot{v}_2s_2}^{-1} \theta_a \pi_{A_1}
\]
is a well defined endomorphism of 
$U_{\Aavv} \cap v_1(U_1)$ and $\til{\sigma}_1^{k+1} \equiv e$
for some sufficiently large integer $k.$ For a given subset
$\Phi \subset \Delta_1^+$ with the property 
\begin{equation}
\lb{prop+}
\alpha, \beta \in \Phi \;\; \mbox{and} \;\; 
\alpha + \beta \in \Delta^+_1 \Rightarrow \alpha + \beta \in \Phi
\end{equation}
define
\[
U^{\Phi}_1= \prod_{\alpha \in \Phi} U^{\alpha}_1.
\]
Recall that \eqref{prop+} implies that 
this product does not depend on the order in which it is taken
and moreover $U^{\Phi}_1$ is an algebraic subgroup of $U_1.$

Extend by linearity $a \colon A_1 \ra A_2$ and $c \colon C_1 \ra C_2$
to bijections $a \colon \De_{A_1}^+ \ra \De_{A_2}^+$ and 
$c \colon \De_{C_1}^+ \ra \De_{C_2}^+.$ For a partial map
$f \colon \De_1 \ra \De_1$ and a subset $Y \subset \Delta_1$
denote $f(Y) := f(Y \cap \Dom (f)),$ where $\Dom (f)$ is the 
domain of $f$. Denote
$d= v_1 c^{-1} v_{2}^{-1} a$ which will be treated as a 
partial bijection from $\De_1$ to $\De_1.$ Set 
\[
\Phi^{(0)}= \left(\Delta^+_1 - \De_{\Aavv}^+ \right)
\cap v_1(\Delta^+_1) 
\;\; \mbox{and inductively} \;\; 
\Phi^{(i+1)}= d(\Phi^{(i)}), \; i \geq 0.
\]  
Since the partial map $d$ preserves addition, all sets
$\Phi^{(i)}$ satisfy the property \eqref{prop+}
and thus
\[
U^{(i)}_1:= U^{\Phi^{(i)}}_1
\]
are algebraic subgroups of $U_1.$ Moreover
$\Phi^{(i+1)} \subset \Phi^{(i)}$ for $i \geq 0$
because $v_1 \in W_1^{C_1}$ and for some positive 
integer $k,$ $\Phi^{(k+1)}= \emptyset$ because 
all powers $d^i$ are well defined on a root
$\alpha \in \Delta_1$ only if 
$\alpha \in \De_{\Aavv}.$ Thus we have the decreasing
sequence of unipotent groups
\begin{equation}
\lb{dec_seq}
U_{\Aavv} \cap v_1(U_1) = U^{(0)}_1 \supset
U^{(1)} \supset \ldots \supset U^{(k)} \supset 
U^{(k+1)}= \{e\}.
\end{equation}     
By direct computation one checks that 
for $\alpha \in \Phi^{(0)},$ 
$\til{\sigma}_1 (U^\alpha_1) = U^{d(\alpha)}_1$
if $d(\alpha)$ is defined 
and $\til{\sigma}_1 (U^\alpha_1) = \{e\}$ otherwise.
As a consequence, 
$\til{\sigma}_1^i(U_{\Aavv} \cap v_1(U_1)) = U^{(i)}_1$
for all $i \geq 1.$ Now \eqref{dec_seq} immediately implies 
that $\til{\sigma}_1$ is a well defined endomorphism of 
$U_{\Aavv} \cap v_1(U_1)$ and $\til{\sigma}_1^{k+1} \equiv e.$

Since $\sigma_1 = \Ad_{m_1} \til{\sigma}_1$ the same will
be true for $\sigma_1$ once we show that 
$M_{\Aavv}$ normalizes all groups $U^{(i)}_1.$ To show this
recall that for a subset $\Phi \subset \Delta^+_1$
satisfying \eqref{prop+}, the unipotent group 
$U^\Phi_1 \subset U_1$ is normalized by $M_{D_1}$ 
(for $D_1 \subset \Ga_1)$ if and only 
\[
\alpha \in \Phi, \; \beta \in \De_{D_1} \; \; 
\mbox{and} \; \; \alpha + \beta \in \Delta^+_1 
\Rightarrow \alpha + \beta \in \Phi.
\]
By induction on $i\geq 0$ one easily checks
that this property is satisfied by $D_1=\Aavv$
and $\Phi= \Phi^{(i)}$ because $\Aavv$ is 
$d$-stable and $\Phi^{(i)}= d(\Phi^{(i-1)}).$ 
\end{proof}

\noindent
To complete the {\em{proof of \thref{main2},}}  
it remains to deduce the formula for $\Stab_U(q).$
Since $v_1 \in {}^{\Aavv}\!W_{1}^{C_1}$ and
$v_2 \in {}^{A_2}\!W_{2}^{\Cbvv}$, it follows from 
\eqref{PP2} and \eqref{PP3} that we have decompositions
\begin{eqnarray}
\lb{U1}
U_{\Aavv} \cap v_1(U_1) & = &\left(U_1 \cap v_1(U_{C_1})\right)\;
 v_1\left(U_{c^{-1}(\Cbvv)}^{C_1}\right)  \\
\lb{U2}
U_{\Cbvv} \cap v_2^{-1}(U_2) &= &\left(U_2 \cap v_2^{-1}(U_{A_2}) 
\right)\;
v_2^{-1}\left(U_{a(\Aavv)}^{A_2}\right).
\end{eqnarray}
Both groups on the right hand side of 
\eqref{U1} are invariant under $\Ad_{m}$ for any $m \in M_{\Aavv}$
because $\Ad_{\dot{v}_1}^{-1} M_{\Aavv} = M_{c^{-1}(\Cbvv)}
\subset M_{C_1}.$ Note also that
\begin{eqnarray}
\lb{phiUC}
&\phi_1(U_{\Cbvv} \cap v_{2}^{-1}(U_2)) \subset v_1
\left(U_{c^{-1}(\Cbvv)}^{C_1}\right) \subset U_{\Aavv} \cap v_1(U_1)\\
\lb{phiUA}
&\phi_2(U_{\Aavv} \cap v_1(U_1)) \subset v_{2}^{-1}
\left(U_{a(\Aavv)}^{A_2}\right) \subset U_{\Cbvv} \cap v_{2}^{-1}(U_2).
\end{eqnarray}
Let now
\[
(n_1u_1, n_2u_2) \in \left(U_{\Aavv} \cap v_1(U_1)\right)
\times \left(U_{\Cbvv} \cap v_{2}^{-1}(U_2)\right)
\]
be decomposed
according to \eqref{U1} and \eqref{U2}.
Then $(n_1u_1, n_2u_2) \in \Stab(q)$ if and only if
\begin{equation}
\lb{un}
\begin{cases}
\Ad_{\dot{v}_2s_2}^{-1}\theta_a\pi_{A_1}(n_1u_1) = u_2&\\
\Ad_{m_1\dot{v}_1}^{-1} (u_1) = \theta_{c}^{-1}\pi_{C_2}
(n_2u_2)&,
\end{cases}
\hs \mbox{or} \hs
\begin{cases}
\phi_2(n_1u_1) = u_2 &\\
u_1 = \phi_1(n_2u_2)&.
\end{cases}
\end{equation}
It follows that 
\begin{equation}
\label{un2}
u_1 =  \phi_1(n_2)\sigma_1(n_1u_1) \hs \mbox{and}
\hs u_2=  \phi_2(n_1)\sigma_2(n_2u_2).
\end{equation}
Let $k$ be such that $\sigma_{i}^{k+1}\equiv e$ for 
$i = 1, 2$. Using the fact that $\sigma_1$ and $\sigma_2$ are group 
homomorphisms and by iterating \eqref{un2}, we get
\begin{equation}
\lb{un3}
u_1 = \psi_1(n_1, n_2)  \hs \mbox{and} \hs
u_2 = \psi_2 (n_1, n_2),
\end{equation}
where $\psi_1$ and $\psi_2$ are respectively given in \eqref{psi1}
and \eqref{psi2}.
Conversely, it is easy to see that any $(n_1u_1, n_2u_2)$ 
of the form in
\eqref{un3} for $(n_1, n_2) \in 
\left(U_1 \cap v_1(U_{C_1}\right)
\times  \left(U_{2} \cap v_{2}^{-1}(U_{A_2})\right)$ satisfies
\eqref{un}. 
This completes the proof of \thref{main2}.
\qed

\bre{normalize}
The Decompositions in \eqref{U1} and \eqref{U2}
are semi-direct products because 
$ v_1\left(U_{c^{-1}(\Cbvv)}^{C_1}\right)$ normalizes 
$U_1 \cap v_1(U_{C_1})$ and 
$v_{2}^{-1}\left(U_{a(\Aavv)}^{A_2}\right)$
 normalizes 
$U_2 \cap v_{2}^{-1}(U_{A_2})$. On the other hand, because
of \eqref{phiUC} and \eqref{phiUA},   the subgroup
\[
\UStab(q) \subset (U_{\Aavv} \cap v_1(U_1)) \times
(U_{\Cbvv} \cap v_2^{-1}(U_2))
\]
can be regarded
as the ``graph" 
\[
\UStab(q) = \{ n \psi(n)\mid \; n \in 
\left( U_1 \cap v_1(U_{C_1}) \right)
\times \left(U_{2} \cap v_{2}^{-1}(U_{A_2})\right)\}
\]
of the map
\[
\psi = \psi_1 \times \psi_2 \colon
\left( U_1 \cap v_1(U_{C_1}) \right)
\times \left(U_{2} \cap v_{2}^{-1}(U_{A_2})\right) \to 
v_1\left(U_{c^{-1}(\Cbvv)}^{C_1}\right) \times 
v_{2}^{-1}\left(U_{a(\Aavv)}^{A_2}\right).
\]
We will show in $\S$\ref{YBE} that the map $\psi$ defines a  
solution to the set-theoretical quantum Yang--Baxter equation
on 
\[
\Stab_U(q) \; \cong \; 
(U_1 \cap v_1(U_{C_1})) 
\times  (U_{2} \cap v_{2}^{-1}(U_{A_2})),
\]
(as algebraic varieties).
\ere

\subsection{Proof of \thref{dim}}\lb{prdim}
By Lemmas \ref{lquintuple} and \ref{lproperty}, 
we know that 
\begin{align*}
\dim R_\A &= \dim U_{A_1} + \dim U_{A_2} + \dim K 
\\
&= \dim P_{A_1} - \dim M_{A_1(v_1,v_2)} + 
\dim U_{A_2} + \dim \Kvv \quad \mbox{and}   
\\
\dim R_\C &= \dim U_{C_1} + \dim U_{C_2} + \dim L 
\\
&=\dim U_{C_1} + \dim P_{C_2} - \dim M_{C_2(v_1,v_2)}
+ \dim \Lvv.
\end{align*}
The second equalities above follow from the fact that 
in the setting of \thref{dim} for any subset $D_1$ of $A_1$
\[
\dim K - \dim K \cap (M_{D_1}\times M_{a(D_1)})
= \dim M_{A_1} - \dim M_{D_1}.
\] 
This is easily deduced from the description of a-graphs
given in \leref{quintuple}.
The second part of \leref{int} and the
part of \thref{main2} on $\UStab(q)$ imply
\begin{align*}
\dim \UStab(q) &= 
\dim U_1 \cap v_1(U_{C_1}) +
\dim U_2\cap v_2^{-1}(U_{A_2})
\\
&= \dim U_{C_1} + \dim U_{A_2} 
-l(v_1) - l(v_2).
\end{align*}
Applying \thref{main2} for the structure 
of $\Stab(q),$ we get 
\begin{align*}
&\dim [m_1 \dot{v}_1, \dot{v}_2 \dot{s}] 
= \dim R_\A + \dim R_\C - \dim \Stab(q) 
\\
&= l(v_1) +l(v_2) + 
\dim P_{A_1} - \dim M_{A_1(v_1,v_2)} + 
\dim P_{C_2} - \dim M_{C_2(v_1,v_2)}
\\
&+\dim \Kvv + \dim \Lvv 
- \dim J_{m_1}(v_1,v_2)
-\dim (\ker \eta_1|_L) \cap (\Id \times v_2^{-1}) (\ker \eta_1|_K).
\end{align*}
Since $\dim \Jvv. m_1 = \dim \Jvv - \dim J_{m_1}(v_1, v_2),$
to complete the proof of \thref{dim} it only remains to show that 
\begin{align}
\lb{final_dim}
&\dim \Jvv = \dim \Kvv + \dim \Lvv - \dim M'_{C_2(v_1,v_2)}
\\
\nn
&-\dim \eta_2 \left( KL \cap 
(v_1^{-1} Z_{A_1(v_1,v_2)} \times v_2 Z_{C_2(v_1,v_2)}) \right) 
-\dim (\ker \eta_1|_L) \cap (\Id \times v_2^{-1}) (\ker \eta_1|_K).
\end{align}
Define 
\[
Q = \{ (k_1, k_2, l_1, l_2) \in \Kvv \times \Lvv \mid
k_2 = l_2 \} \subset \Kvv \times \Lvv.
\]
In the two cases of \thref{dim} $Q$ is respectively an algebraic/Lie 
subgroup of $\Kvv \times \Lvv.$ By definition $\Jvv$ is the
image of $Q$ under $(k_1, k_2, l_1, l_2) \mapsto (k_1, l_1),$
in particular $\Jvv$ is also an algebraic/Lie group. It is easy to see
that the kernel of this homomorphism is 
$(\ker \eta_1|_L) \cap (\Id \times v_2^{-1}) (\ker \eta_1|_K).$ 
Therefore
\begin{equation}
\lb{final_dim1}
\dim \Jvv = \dim Q - \dim (\ker \eta_1|_L) \cap 
(\Id \times v_2^{-1}) (\ker \eta_1|_K).
\end{equation} 
The homogeneous space $(\Kvv \times \Lvv)/Q$ is isomorphic
to $\eta_2(\Kvv)\eta_2(\Lvv) \subset M_{C_2(v_1,v_2)}$
by $(k_1, k_2, l_1, l_2)Q \mapsto k_2 l_2^{-1}.$ 
Moreover from \leref{quintuple} one obtains 
\[
\eta_2(\Kvv)\eta_2(\Lvv) = M'_{C_2(v_1,v_2)}
\eta_2\left(\Kvv \Lvv \cap (Z_{A_1(v_1,v_2)} 
\times Z_{C_2(v_1,v_2)}) \right) 
\]
and
\begin{align}
\lb{final_dim2}
\dim Q &= \dim \Kvv + \dim \Lvv - \dim M'_{C_2(v_1,v_2)} 
\\
\nn
&- \dim \eta_2\left(\Kvv \Lvv \cap (Z_{A_1(v_1,v_2)}
\times Z_{C_2(v_1,v_2)}) \right).
\end{align}
Finally substituting \eqref{final_dim1} in \eqref{final_dim2} 
and using that the projections of 
$\Kvv \Lvv \cap (Z_{A_1(v_1,v_2)}
\times Z_{C_2(v_1,v_2)})$ and
$K L \cap (v_1^{-1}Z_{A_1(v_1,v_2)} \times v_2 Z_{C_2(v_1,v_2)})$ 
under $\eta_2$ have the same dimensions lead to \eqref{final_dim}.
\qed
\sectionnew{Solutions to the set-theoretical 
quantum Yang--Baxter equation}
\lb{YBE}

Recall that for a set $V$, the set-theoretical 
quantum Yang--Baxter equation
for an invertible map $T \colon V \times V \to V \times V$ is 
\begin{equation}
\lb{YB}
T^{12}T^{13}T^{23} = T^{23}T^{13}T^{12},
\end{equation}
where for $i, j \in \{1, 2, 3\}$, $T^{ij}$ 
denotes the map from $V \times V \times V$ 
to itself that has $T$ act on the
$(i,j)$'th components. 
 
 For $v_1 \in W_{1}^{C_1}$ and $v_2 \in {}^{A_2}\!W_{2}$,
we will set, for notational simplicity,
\begin{align*}
N_{v_1, v_2} &=\left( U_1 \cap v_1(U_{C_1}) \right)
\times \left(U_{2} \cap v_{2}^{-1}(U_{A_2})\right)\\
Q_{v_1, v_2} & = v_1\left(U_{c^{-1}(\Cbvv)}^{C_1}\right) \times 
v_{2}^{-1}\left(U_{a(\Aavv)}^{A_2}\right).
\end{align*}
Then $Q_{v_1, v_2}$ normalizes $N_{v_1, v_2}$, and
\[
N_{v_1, v_2} Q_{v_1, v_2} = 
(U_{\Aavv} \cap v_1(U_1)) \times
(U_{\Cbvv} \cap v_2^{-1}(U_2))
\]
Recall from \reref{normalize}
that for a point
$q = (\dot{v}_1m, \dot{v}_2 s_2) \in G \times G$, where,
$\dot{v}_i$ is a representative of $v_i$ in the normalizer
of $H_i$ in $G_i$ for $i = 1, 2$, $m \in M_{\Aavv}$, and
$s_2 \in Z_{\Cbvv}$, we have
\[
\Stab_U(q) = \{n \psi(n) \mid n \in N_{v_1, v_2}\} \subset
(U_{\Aavv} \cap v_1(U_1)) \times
(U_{\Cbvv} \cap v_2^{-1}(U_2)),
\]
where $\psi = \psi_1 \times \psi_2 \colon N_{v_1, v_2} \to  
Q_{v_1, v_2}$
and  $\psi_1$ and $\psi_2$ are respectively given in
\eqref{psi1} and \eqref{psi2}.
In this section, we show that the map $\psi$ defines a 
solution to the 
set-theoretical quantum Yang--Baxter equation on $N_{v_1, v_2}$. 

As is true for any group \cite{D-set}, the map
\begin{equation}
\lb{T0}
T_0 \colon N_{v_1, v_2} \times N_{v_1, v_2} \to 
N_{v_1, v_2} \times N_{v_1, v_2} \colon 
(n, \; n^\prime) \mapsto (n^\prime, \;  (n^\prime)^{-1} n n^{\prime})
\end{equation}
satisfies the set-theoretical quantum Yang--Baxter equation
\eqref{YB}.
We now show that a twist of $T_0$ using the map $\psi$ is also
a solution to \eqref{YB}.

Let $\sigma$ and $F$ be the maps from $N_{v_1, v_2} \times N_{v_1, v_2}$
to itself given respectively by
\begin{align}
\lb{tau}
&\sigma(n, n^\prime) = (n^\prime, n),\\
\lb{F}
&F(n, n^\prime) = (n, \; \psi(n) n^\prime \psi(n)^{-1}).
\end{align}
Define 
\begin{equation}
\lb{T}
T = (\sigma F \sigma)^{-1} T_0 \,F \colon N_{v_1, v_2} \times N_{v_1, v_2}
\to N_{v_1, v_2} \times N_{v_1, v_2}.
\end{equation}
 
\bpr{YB}
The map $T$ is 
a solution to the set-theoretical quantum Yang--Baxter equation 
\eqref{YB} on $N_{v_1, v_2}$.
\epr
 
\prref{YB} is a special case of the following general fact.

\ble{YBE}
Assume that $U$ is a group and that $Q$ and $N$ are subgroups of $U$ 
such that
$N \cap Q = \{e\}$ and that $Q$ normalizes $N$. Suppose that
$\psi \colon N \to Q$ is a map such that  
\[
S := \{ n\psi(n)\mid n \in N\} 
\]
is a subgroup of $U$. Let $T_0$, $\sigma$, $F$, 
and $T$ be the maps from $N \times N$ 
to itself given  by
replacing $N_{v_1, v_2}$ by $N$ in \eqref{T0}, \eqref{tau}, \eqref{F},
and \eqref{T} respectively. Then  $T$ is a 
solution to the set-theoretical quantum Yang--Baxter Equation 
\eqref{YB} on $N$.
\ele

\begin{proof} Since $S$ is a subgroup of $U$, the map
\[
\wt{\psi} \colon S \rightarrow Q \colon n\psi(n) \mapsto \psi(n)
\]
is a group homomorphism. Thus the map
\[
S \times N \rightarrow N \colon (n\psi(n), n^\prime) \mapsto
(n\psi(n)) \cdot n^\prime :=\psi(n) n^\prime \psi(n)^{-1}
\]
defines a left action of $S$ on $N$ by group automorphisms. 
Moreover, it is clear from the definition that
the map $\pi \colon S \to N \colon n\psi(n) \mapsto n$ is a 
bijective $1$-cocycle on $S$
with coefficients in $N$, i.e., $\pi$ is bijective and that 
\[
\pi(s_1 s_2) = \pi(s_1) (s_1 \cdot \pi(s_2)), \hspace{.2in}
s_1, s_2 \in S.
\]
By \cite[Theorem 6]{LYZ}, the map $T$ is a 
solution to \eqref{YB}.
\end{proof}

\bre{trivial-twist}
Although for any integer $m \geq 2$, the action of the Braid group $B_m$ 
on $N^{\times m}$
induced by $T$ is isomorphic to the one induced by $T_0$ 
\cite[Theorem 6]{LYZ}, 
the fact that
$T$, as the twisting of $T_0$ by $F$, still satisfies the 
quantum Yang--Baxter
Equation \eqref{YB} is non-trivial.
When $N$ is abelian, the fact that $T$ satisfies \eqref{YB}
is also proved in \cite{ESSo}. 
\ere
\sectionnew{The groups $R_\A$ are spherical subgroups of $G_1\times G_2$} 
\lb{spherical}

Recall that an algebraic subgroup $R$ of a reductive group $G$ is called 
spherical if $R$ has finitely many orbits on the 
flag variety $G/B$ where $B$ is a Borel subgroup of $G.$ 
In this section, we fix an admissible pair 
$\A = (a, K)$ for $G_1 \times G_2$ where $K$ is
an algebraic subgroup of $G_1 \times G_2.$
Let $A_1$ and $A_2$ be the domain and the range of $a$ respectively.
For $i = 1, 2$, we fix a representative $\dot{v}_i$ in the normalizer
of $H_i$ in $G_i$ for each $v_i \in W_i$.
The main result in this section is:

\bpr{RAPP}  
Let $C_1 \subset \Ga_1$ and $C_2 \subset \Ga_2$
be arbitrary.
Then every $(R_\A, P_{C_1} \times P_{C_2})$ double 
coset in $G_1 \times G_2$ contains
exactly one point of the form $(\dot{v}_1, \dot{v}_2)$, where
$v_1 \in {}^{A_1}\!W_{1}^{C_1},$ 
$v_2 \in {}^{a(A_1 \cap v_1(C_1))}\!W_{2}^{C_2}$.
\epr

As an immediate corollary of \prref{RAPP} we obtain:

\bco{spher} 
If $(a, K)$ is an admissible pair for $G_1 \times G_2$ 
such that $K$ is an algebraic subgroup
of $G_1 \times G_2,$ then the 
group $R_\A$ is a spherical subgroup of $G_1 \times G_2.$
\eco 

We will first prove an auxiliary Lemma which can be also viewed as
a special case of \prref{RAPP}. 
Recall that for $D_i \subset A_i \subset \Ga_i$, 
$P_{D_i}^{A_i} = P_{D_i} \cap M_{A_i}$ denotes 
the standard parabolic subgroup of $M_{A_i}$ determined by $D_i$.

\ble{special-1}
For any subsets $D_1 \subset A_1$ and 
$D_2 \subset A_2$, 
every $(K, P_{D_1}^{A_1} \times P_{D_2}^{A_2})$ double coset
in $M_{A_1} \times M_{A_2}$ contains exactly one point of the form 
$(e,\dot{w})$, where $w \in {}^{a(D_1)}\!W_{A_2}^{D_2}$.
\ele

\begin{proof} Let $H_{A_i} = H_i \cap M_{A_i}$ for $i = 1, 2$.
Since $\eta_1(K) \supset M'_{A_1}$ and
$P_{D_1}^{A_1} \supset H_{A_1}$, 
every $(K,  P_{D_1}^{A_1} \times P_{D_2}^{A_2})$ double coset
in $M_{A_1} \times M_{A_2}$
is of the form $[e, m]$ for some $m \in M_{A_2}$. 
Let $\tP_{a(D_1)}^{A_2} = P_{a(D_1)}  \cap \eta_2(K)$. 
Since $H_{A_2} \tP_{a(D_1)}^{A_2} = P_{a(D_1)}^{A_2}$ 
the  Bruhat decomposition of $M_{A_2},$
see \eqref{eq_G-Bruhat-Pa}, implies
\[
M_{A_2} = \coprod_{w \in {}^{a(D_1)}W_{A_2}^{D_2}}
\tP_{a(D_1)}^{A_2} \dot{w} P_{D_2}^{A_2}.
\]
For $m \in M_{A_2}$, write $m = q_2 \dot{w} p_2$ for 
$w \in {}^{a(D_1)}W_{A_2}^{D_2}$,
$q_2 \in \tP_{a(D_1)}^{A_2}$, and 
$p_2 \in P_{D_2}^{A_2}$. Then by Part 3) of
\leref{property}, we know that there exists $q_1 \in P_{D_1}^{A_1}$
such that $(q_1, q_2) \in K$. Thus
\[
[e, m] = [e, q_2\dot{w}p_2]=
[q_{1}^{-1}, \dot{w}p_2] = [e, \dot{w}].
\]
It is easy to see that if $w, w' \in {}^{a(D_1)}W_{A_2}^{D_2}$ 
are such that
$[e, \dot{w}]=[e,\dot{w}']$, then $w = w'$ 
which completes the proof of the lemma.
\end{proof}

\noindent
{\em{Proof of \prref{RAPP}.}}
For the left action of  $P_{A_1} \times P_{A_2}$ on
$(G_1 \times G_2)/(P_{C_1} \times P_{C_2})$, we know from
the Bruhat decomposition
that every $P_{A_1} \times P_{A_2}$ orbit passes through a point
$(\dot{w}_1, \dot{w}_2)(P_{C_1} \times P_{C_2})$ for a unique
$w_1 \in {}^{A_1}\!W^{C_1}$ and a unique $w_2 \in {}^{A_2}\!W^{C_2}$.
Denote by $\Stab_{(\dot{w}_1, \dot{w}_2)}$ the stabilizer subgroup of
$P_{A_1} \times P_{A_2}$ at the point
$(\dot{w}_1, \dot{w}_2)(P_{C_1} \times P_{C_2})$.
Then by \leref{abs_ind}, every $R_\A$-orbit in
$(G_1 \times G_2)/(P_{C_1} \times P_{C_2})$ contains a point of the form
$(p_1 \dot{w}_1, p_2 \dot{w}_2)(P_{C_1} \times P_{C_2})$ for a unique
$w_1 \in {}^{A_1}\!W^{C_1}$, a unique $w_2 \in {}^{A_2}\!W^{C_2}$, and
for some $(p_1, p_2) \in P_{A_1} \times P_{A_2}$, and two such points
$(p_1 \dot{w}_1, p_2 \dot{w}_2)(P_{C_1} \times P_{C_2})$ and
$(q_1 \dot{w}_1, q_2 \dot{w}_2)(P_{C_1} \times P_{C_2})$  are in the same
$R_\A$-orbit
if and only if
$(p_1, p_2)$ and $(q_1, q_2)$ are in the same
$R_\A$-orbit on $(P_{A_1} \times P_{A_2})/\Stab_{(\dot{w}_1, \dot{w}_2)}$.
Thus we need to understand the double coset spaces
$R_\A \backslash (P_{A_1} \times P_{A_2})/\Stab_{(\dot{w}_1, \dot{w}_2)}$.

Fix now $w_1 \in {}^{A_1}\!W^{C_1}$ and  $w_2 \in {}^{A_2}\!W^{C_2}$.
Then it is easy to see that
\[
\Stab_{(\dot{w}_1, \dot{w}_2)} = (P_{A_1} \cap w_1(P_{C_1})) \times
(P_{A_2} \cap w_2(P_{C_2})).
\]
On the other hand, recall that for $i = 1, 2$, 
$\pi_{A_i} \colon P_{A_i} \to M_{A_i}$ is 
the projection with respect to the
decomposition $P_{A_i} = M_{A_i}  N_{A_i}$. The projection
$\pi_{A_1} \times \pi_{A_2} \colon P_{A_1} \times P_{A_2} \to
M_{A_1} \times M_{A_2}$ induces an isomorphism between the
quotient
spaces $R_\A \backslash (P_{A_1} \times P_{A_2})$ and
$K \backslash (M_{A_1} \times M_{A_2})$. We know from
1) of \leref{int} that
\[
(\pi_{A_1} \times \pi_{A_2})(\Stab_{(\dot{w}_1, \dot{w}_2)}) =
P_{A_1 \cap w_1(C_1)}^{A_1} \times P_{A_2 \cap w_2(C_2)}^{A_2}.
\]
Since $M_{A_1} \times M_{A_2}$ normalizes  $N_{A_1} \times N_{A_2}$,
we see that
the projection
$\pi_{A_1} \times \pi_{A_2}$ induces an isomorphism of
double coset spaces
\[
R_\A \backslash P_{A_1} \times P_{A_2} / \Stab_{(\dot{w}_1, \dot{w}_2)}
\; \cong \; K \backslash  M_{A_1} \times M_{A_2}/
(P_{A_1 \cap w_1(C_1)}^{A_1} \times P_{A_2 \cap w_2(C_2)}^{A_2}).
\]
By \leref{special-1}, every
$(K, P_{A_1 \cap w_1(C_1)}^{A_1} \times P_{A_2 \cap
w_2(C_2)}^{A_2})$-double
coset in $M_{A_1} \times M_{A_2}$ contains exactly one point
of the form $(e, \dot{w})$ for some
$w \in {}^{a(A_1 \cap w_1(C_1))}W_{A_2}^{A_2 \cap w_2(C_2)}$.
Thus every $R_\A$-orbit in $(G_1 \times G_2)/( P_{C_1} \times P_{C_2})$
contains a unique element of the form $(\dot{w}_1, \dot{w}\dot{w}_2)$, 
where $w_1 \in  {}^{A_1}\!W^{C_1}$, $w_2 \in {}^{A_2}\!W^{C_2}$, and
$w \in {}^{a(A_1 \cap w_1(C_1))}W_{A_2}^{A_2 \cap w_2(C_2)}$.
Let $v_1 = w_2$ and $v_2 = ww_2$. \prref{RAPP} now follows from
\prref{W_prod}.
\qed

\end{document}